\documentclass[12pt]{article}
\usepackage{amsmath,amsfonts,amssymb,amsthm}
\input amssym.def
\begin{document}
\def \Z{\Bbb Z}
\def \C{\Bbb C}
\def \R{\Bbb R}
\def \Q{\Bbb Q}
\def \N{\Bbb N}

\def \A{{\mathcal{A}}}
\def \D{{\mathcal{D}}}
\def \E{{\mathcal{E}}}
\def \H{\mathcal{H}}
\def \S{{\mathcal{S}}}
\def \wt{{\rm wt}}
\def \tr{{\rm tr}}
\def \span{{\rm span}}
\def \Res{{\rm Res}}
\def \Der{{\rm Der}}
\def \End{{\rm End}}
\def \Ind {{\rm Ind}}
\def \Irr {{\rm Irr}}
\def \Aut{{\rm Aut}}
\def \GL{{\rm GL}}
\def \Hom{{\rm Hom}}
\def \mod{{\rm mod}}
\def \ann{{\rm Ann}}
\def \ad{{\rm ad}}
\def \rank{{\rm rank}\;}
\def \<{\langle}
\def \>{\rangle}

\def \g{{\frak{g}}}
\def \h{{\hbar}}
\def \k{{\frak{k}}}
\def \sl{{\frak{sl}}}
\def \gl{{\frak{gl}}}

\def \be{\begin{equation}\label}
\def \ee{\end{equation}}
\def \bex{\begin{example}\label}
\def \eex{\end{example}}
\def \bl{\begin{lem}\label}
\def \el{\end{lem}}
\def \bt{\begin{thm}\label}
\def \et{\end{thm}}
\def \bp{\begin{prop}\label}
\def \ep{\end{prop}}
\def \br{\begin{rem}\label}
\def \er{\end{rem}}
\def \bc{\begin{coro}\label}
\def \ec{\end{coro}}
\def \bd{\begin{de}\label}
\def \ed{\end{de}}

\newcommand{\m}{\bf m}
\newcommand{\n}{\bf n}
\newcommand{\nno}{\nonumber}
\newcommand{\nord}{\mbox{\scriptsize ${\circ\atop\circ}$}}
\newtheorem{thm}{Theorem}[section]
\newtheorem{prop}[thm]{Proposition}
\newtheorem{coro}[thm]{Corollary}
\newtheorem{conj}[thm]{Conjecture}
\newtheorem{example}[thm]{Example}
\newtheorem{lem}[thm]{Lemma}
\newtheorem{rem}[thm]{Remark}
\newtheorem{de}[thm]{Definition}
\newtheorem{hy}[thm]{Hypothesis}
\makeatletter \@addtoreset{equation}{section}
\def\theequation{\thesection.\arabic{equation}}
\makeatother \makeatletter

\begin{center}
{\Large \bf On a category of $\gl_{\infty}$-modules}
\end{center}

\begin{center}
{Cuipo Jiang\footnote{Supported by China NSF grant 11371245,  China RFDP grant 2010007310052, and the Innovation Program of Shanghai Municipal Education Commission (11ZZ18)}\\
Department of Mathematics\\
Shanghai Jiaotong University, Shanghai 200240, China\\
Haisheng Li\footnote{Partially supported by NSA grant
H98230-11-1-0161 and China NSF grant 11128103}\\
Department of Mathematical Sciences\\
Rutgers University, Camden, NJ 08102, USA}
\end{center}

\begin{abstract}
We study a particular category ${\cal{C}}$ of $\gl_{\infty}$-modules and a
subcategory ${\cal{C}}_{int}$ of integrable $\gl_{\infty}$-modules.
As the main results, we classify the irreducible modules in these
two categories and we show that every module in
category ${\cal{C}}_{int}$ is semi-simple. Furthermore, we determine
the decomposition of the tensor products of irreducible modules in
category ${\cal{C}}_{int}$.
\end{abstract}

\section{Introduction}
In the representation theory of Kac-Moody algebras, of great importance is
the category of integrable modules, where classifying 
irreducible integrable modules has been an open problem (see \cite{kac1}). 
For (finite rank) affine Kac-Moody algebras,
irreducible integrable modules with finite-dimensional
weight spaces were classified by Chari (see \cite{ch}).
Integrable modules for affine Kac-Moody algebras were studied further
in \cite{cp1,cp2,cp3} (cf. \cite{rao1}, \cite{li-mz}).

Integrable representations for infinite rank affine Kac-Moody algebras, including
$\gl_{\infty}$ (the Lie algebra of
doubly infinite matrices with only finitely many nonzero entries),
are also of great importance in many areas, especially in mathematical physics.
Among the most important and interesting 
results is the remarkable relation of highest weight integrable representations with
soliton equations, which was discovered and developed by Kyoto
school (see \cite{dkm}, \cite{djkm1, djkm2}, \cite{jm}). 
Lie algebra $\gl_{\infty}$ has also been used effectively
to study ${\mathcal{W}}$-algebras (see \cite{fkrw}, \cite{kr}).

In a previous study \cite{jli}, we exhibited a natural association
of quantum vertex algebras (see \cite{li-qva1}) to $\gl_{\infty}$. In
that study we came across a category of $\gl_{\infty}$-modules $W$
satisfying the condition that for any $m\in \Z,\ w\in W$,
$E_{m,n}w=0$ for all but finitely many integers $n$.  (A canonical
base of $\gl_{\infty}$ consists of $E_{m,n}$ $(m,n\in \Z)$, where
$E_{m,n}$ denotes the matrix whose only nonzero entry is the
$(m,n)$-entry which is $1$.) This category contains the natural
module $\C^{\infty}$ and its tensor
products, whereas it excludes nontrivial highest
weight modules and lowest weight modules. It is our hope to classify
the irreducible objects in this category. This is the main motivation for this current paper.

In this paper,  we focus on a smaller category of
$\gl_{\infty}$-modules, for which we are able to determine and
classify all the irreducible objects. Specifically, we study
$\gl_{\infty}$-modules $W$ satisfying the condition that for every
$w\in W$, there exists a finite subset $S$ of $\Z$ such that
$$E_{m,n}w=0\ \ \ \mbox{ for }m,n\in \Z,\ n\notin S.$$
Denote by ${\mathcal{C}}$
the category of such $\gl_{\infty}$-modules and by ${\mathcal{C}}_{int}$
the category of those integrable $\gl_{\infty}$-modules.
The category ${\mathcal{C}}$ still contains the natural
module $\C^{\infty}$, and it is
closed under tensor product.

Note that for a finite-dimensional simple Lie algebra, or more generally for a
Kac-Moody algebra (see \cite{kac1}), one has the well known category $o$ and its subcategory $o_{int}$ of
integrable modules. To a certain extent, this category ${\cal{C}}$
of $\gl_{\infty}$-modules is analogous to category $o$.
Especially, it is proved that every module in category ${\mathcal{C}}_{int}$ is
completely reducible. By using certain generalized
Verma modules, we can classify irreducible modules in
categories ${\mathcal{C}}$ and ${\mathcal{C}}_{int}$.
On the other hand, we can also determine the decomposition of tensor
products of irreducible modules in category ${\mathcal{C}}_{int}$ 
in terms of the decomposition of tensor
products of irreducible modules for $\gl_{n}$ with sufficiently large $n$.

We now give a more detailed description of the main results.
Let $S$ be a finite and nonempty subset
of $\Z$. To $S$ we associate a triangular decomposition
$$\gl_{\infty}=\gl_{\infty}^{(S,+)}\oplus \gl_{\infty}^{(S,0)}\oplus \gl_{\infty}^{(S,-)},$$
where
\begin{eqnarray*}
&&\gl_{\infty}^{(S,+)}={\rm span}\{ E_{m,p}\ |\ m\in \Z,\ p\notin
S\},\ \ \ \ \gl_{\infty}^{(S,-)}= {\rm span}\{ E_{p,n}\ |\ p\notin
S,\ n\in S\},\\
&&\hspace{3cm} \gl_{\infty}^{(S,0)}={\rm span}\{ E_{m,n}\ |\ m,n\in
S\}.
\end{eqnarray*}
Alternatively, set
$$\gl_{S}=\gl_{\infty}^{(S,0)}={\rm span}\{ E_{m,n}\ |\ m,n\in S\}.$$ Given a
$\gl_{S}$-module $U$, using this particular triangular
decomposition we define a generalized Verma module $M(S,U)$,
which is a $\gl_{\infty}$-module induced from $\gl_{S}$-module $U$.
When $U$ is irreducible, $M(S,U)$ has a unique irreducible quotient module,
denoted by $L(S,U)$.
 We show that $M(S,U)$ and
$L(S,U)$ belong to category ${\cal{C}}$ and that every irreducible
$\gl_{\infty}$-module in ${\cal{C}}$ is isomorphic to $L(S,U)$ for
some finite subset $S$ of $\Z$ and for some irreducible
$\gl_{S}$-module $U$. We also give a necessary and sufficient condition that $L(S_{1},U_{1})\simeq L(S_{2},U_{2})$, where $S_{1},S_{2}$ are finite subsets of $\Z$
and $U_{1},U_{2}$ are irreducible modules for $\gl_{S_{1}}$ and $\gl_{S_{2}}$, respectively.

Set
$$H={\rm span}\{ E_{n,n}\ |\ n\in \Z\},$$ a Cartan subalgebra of
$\gl_{\infty}$.
Let $S$ be a finite nonempty subset of $\Z$. Set
$$H_{S}={\rm span}\{ E_{n,n}\ |\ n\in S\},$$ which is a Cartan subalgebra of
$\gl_{S}$.
Let $\lambda\in H^{*}$ such that
$\lambda(E_{n,n})=0$ for $n\notin S$.
Define $M(S,\lambda)$ to be
the generalized Verma module $M(S,U)$ with $U=M(\lambda_{S})$,
where $M(\lambda_{S})$ denotes the Verma $\gl_{S}$-module with highest weight $\lambda_{S}=\lambda_{H_{S}}$.
A fact is that the weight-$\lambda$ subspace of $M(S,\lambda)$ is
one-dimensional. We then define an irreducible module
$L(S,\lambda)$ as the quotient module of
$M(S,\lambda)$ by the maximal submodule. We show that $L(S,\lambda)$ is integrable if and
only if $\lambda(E_{n,n})\in \N$ for $n\in \Z$ and
$\lambda(E_{m,m})\ge \lambda(E_{n,n})$ for $m,n\in S$ with $m<n$.
Furthermore, we show that every such irreducible integrable module
belongs to category ${\cal{C}}_{int}$ and every irreducible
$\gl_{\infty}$-module in ${\cal{C}}_{int}$ is isomorphic to such an
irreducible module.

We furthermore study the tensor products of irreducible modules in ${\mathcal{C}}_{int}$.
It is shown that the tensor product of any two irreducible modules in ${\mathcal{C}}_{int}$
is always cyclic. Let $W$ be any $\gl_{\infty}$-module in ${\mathcal{C}}$
and let $S$ be a finite subset of $\Z$. Set
$$\Omega_{S}(W)=\{ w\in W\ |\ E_{p,q}w=0\ \ \mbox{ for any }p,q\in \Z\ \mbox{with }q\notin S\},$$
which is a $\gl_{S}$-submodule of $W$.
It is proved that if $W$ is irreducible, there exists a finite subset $S'$ of $\Z$ such that $\Omega_{S}(W)$ is an irreducible $\gl_{S}$-module for any finite subset $S$ of $\Z$, containing $S'$. Now, let $W_{1},W_{2}$ be irreducible $\gl_{\infty}$-modules in ${\mathcal{C}}_{int}$.
It is proved that for any sufficiently large $S$, the decomposition of $W_{1}\otimes W_{2}$ into irreducible $\gl_{\infty}$-modules is determined by the decomposition of $\Omega_{S}(W_{1})\otimes \Omega_{S}(W_{2})$ into irreducible $\gl_{S}$-modules. Consequently, the tensor product of any two
irreducible modules in category ${\mathcal{C}}_{int}$ always decomposes into a finite
sum of irreducible modules, unlike the case with highest weight
modules.

After this paper was completed, we found a very interesting paper
\cite{ps2000}, in which Penkov and Serganova had studied the
category of integrable modules with finite dimensional weight
subspaces for $sl(\infty)$, $o(\infty)$, and $sp(\infty)$. Among
other results they proved that the category of integrable modules
with finite dimensional weight subspaces is semisimple and they also
identified each irreducible module in this category. Their method is somewhat different
from that of this present paper.

This paper is organized as follows: In Section 2, we define
categories ${\cal{C}}$ and ${\cal{C}}_{int}$ of
$\gl_{\infty}$-modules and we establish the complete reducibility.
In Section 3, we classify irreducible modules in categories ${\mathcal{C}}$ and ${\mathcal{C}}_{int}$. In Section 4, we study the decomposition of the
tensor product modules.

\section{Categories ${\mathcal{C}}$ and ${\cal{C}}_{int}$ of $\gl_{\infty}$-modules}
In this section, we introduce a category ${\cal{C}}$ of
$\gl_{\infty}$-modules and a subcategory ${\cal{C}}_{int}$ of
integrable modules. As
the main result of this section, we prove
that every $\gl_{\infty}$-module in category ${\cal{C}}_{int}$ is
completely reducible.

We begin with Lie algebra $\gl_{\infty}$, which is the Lie algebra
of doubly infinite matrices with only finitely many nonzero entries,
under the commutator bracket. A canonical base consists of $E_{m,n}$
$(m,n\in \Z)$, where $E_{m,n}$ denotes the matrix whose only nonzero
entry is the $(m,n)$-entry which is $1$, and we have
\begin{eqnarray}
[E_{m,n},E_{r,s}]=\delta_{n,r}E_{m,s}-\delta_{m,s}E_{r,n}
\end{eqnarray}
for $m,n,r,s\in \Z$.
Let $\C^{\infty}$ denote the vector space of doubly infinite column
vectors with only finitely many nonzero entries. Denote the standard
unit base vectors by $v_{n}$ for $n\in \Z$. The natural action of
$\gl_{\infty}$ on $\C^{\infty}$ is given by
$$E_{i,j}v_{k}=\delta_{j,k}v_{i}\ \ \ \mbox{ for }i,j, k\in \Z.$$

Define $\deg E_{i,j}=j-i$ for $i,j\in \Z$ to make $\gl_{\infty}$ a
$\Z$-graded Lie algebra, where the degree-$n$ homogeneous subspace
$(\gl_{\infty})_{(n)}$ for $n\in \Z$ is linearly spanned by
$E_{m,m+n}$ for $m\in \Z$. We have the standard triangular decomposition
$$\gl_{\infty}=\gl_{\infty}^{+}\oplus \gl_{\infty}^{0}\oplus
\gl_{\infty}^{-},$$ where $\gl_{\infty}^{\pm}=\sum_{\pm(j-i)>0}\C
E_{i,j}$ and $\gl_{\infty}^{0}=\sum_{n\in \Z}\C E_{n,n}$. Alternatively, set
\begin{eqnarray}
H=\gl_{\infty}^{0}={\rm span}\{ E_{n,n}\ |\ n\in \Z\},
\end{eqnarray}
a Cartan subalgebra of $\gl_{\infty}$.

The formal completion $\overline{\gl_{\infty}}$ of $\gl_{\infty}$ is
also a $\Z$-graded Lie algebra, where
$$\overline{(\gl_{\infty})}_{(n)}=\left\{ \sum_{m\in \Z} a_{m}E_{m,m+n}\ |\; a_{m}\in \C\right\}.$$
If $W$ is a $\gl_{\infty}$-module such that for every $w\in W$ and
for every $n\in \Z$, $E_{m,m+n}w=0$ for all but finitely many
integers $m$, then $W$ is naturally a
$\overline{\gl_{\infty}}$-module.

\bd{dintegrable} {\em A $\gl_{\infty}$-module $W$ is said to be {\em
integrable} if for every $n\in \Z$, $E_{n,n}$ is semi-simple on $W$
and if for any $p,q\in \Z$ with $p\ne q$, $E_{p,q}$ is locally
nilpotent on $W$.} \ed

A notion of integrable module for a general Lie algebra including
Kac-Moody Lie algebras was introduced in \cite{kac1}.
This definition of an integrable module for $\gl_{\infty}$ is just a
version of that. {}From \cite{kac1} (Proposition 3.8), a
$\gl_{\infty}$-module $W$ is integrable if $E_{n,n}$ is semi-simple
on $W$ for {\em some} $n\in \Z$ and if $E_{j,j+1}$ and $E_{j+1,j}$
for all $j\in \Z$ are locally nilpotent on $W$.

\br{rfd-integrable} {\em Let $\g$ be a finite-dimensional simple Lie
algebra with Chevalley generators $e_{i},f_{i},h_{i}$ $(1\le i\le
l)$. A $\g$-module $W$ is {\em integrable} if $e_{i},f_{i}$ $(1\le
i\le l)$ are locally nilpotent on $W$. From \cite{kac1} (Proposition
3.8), every integrable $\g$-module is $\g$-locally finite, hence a
direct sum of finite-dimensional irreducible modules. Consequently,
on an integrable $\g$-module, every root vector of $\g$ is locally
nilpotent and the Cartan algebra is semi-simple. } \er

\br{rsl2} {\em Consider the three dimensional simple Lie algebra
$\g=\sl(2,\C)$ with the standard Chevalley generators $e,f,h$.
Suppose that $V$ is an integrable $\sl(2,\C)$-module with a nonzero
vector $v$ satisfying the condition that $ev=0$ and $hv=kv$ for some
$k\in \C$. Then $k\in \N$, $f^{k+1}v=0$, and
$e^{k}f^{k}v=(k!)^{2}v$. } \er

Recall that a highest weight $\gl_{\infty}$-module with
highest weight $\lambda\in H^{*}$ is a module $W$
with a vector $w$ such that
\begin{eqnarray*}
&&E_{n,n}w=\lambda_{n}w,\ \
E_{n,n+1}w=0\ \ \ \mbox{ for }n\in \Z,\\
&&W=U(\gl_{\infty})w,
\end{eqnarray*}
where $\lambda_{n}=\lambda(E_{n,n})$. {}From \cite{kac1}, a highest
weight integrable $\gl_{\infty}$-module is irreducible.

\bd{dcategoryD} {\em Denote by ${\cal{C}}$ the category of
$\gl_{\infty}$-modules $W$ such that for any $w\in W$, there exists
a finite subset $S$ of $\Z$ such that $E_{m,n}w=0$ for all $m,n\in
\Z$ with $n\notin S$. Furthermore, define ${\cal{C}}_{int}$ to be
the subcategory consisting of integrable $\gl_{\infty}$-modules in
${\cal{C}}$.} \ed

It can be readily seen that every submodule
of a $\gl_{\infty}$-module from ${\cal{C}}$ is in ${\cal{C}}$ and
the tensor product of any (finitely many) $\gl_{\infty}$-modules
from ${\cal{C}}$ is in ${\cal{C}}$. The same can be said for
category ${\cal{C}}_{int}$.

Set
\begin{eqnarray}
I_{\infty}=\sum_{n\in \Z}E_{n,n},
\end{eqnarray}
which lies in the completion $\overline{\gl_{\infty}}$ of
$\gl_{\infty}$. Notice that every $\gl_{\infty}$-module in category
${\cal{C}}$ is naturally a $\overline{\gl_{\infty}}$-module. Then
for any $\gl_{\infty}$-module $W$ in category ${\cal{C}}$,
$I_{\infty}$ is a well defined operator on $W$ and
$[I_{\infty},\gl_{\infty}]=0$.

Let $S$ be a finite and nonempty subset of $\Z$. Set
\begin{eqnarray}
\gl_{S}={\rm span}\{ E_{m,n}\ |\ m,n\in S\},
\end{eqnarray}
which is a subalgebra of $\gl_{\infty}$, and set
\begin{eqnarray}
I_{S}=\sum_{n\in S}E_{n,n}\in \gl_{S}.
\end{eqnarray}
The Lie algebra $\gl_{S}$ is reductive with one-dimensional center $\C
I_{S}$.

\bl{lenn-eigenvalues} Let $W$ be any $\gl_{\infty}$-module in
category ${\cal{C}}_{int}$. Then for every $n\in \Z$, $E_{n,n}$ is
semi-simple on $W$ with only nonnegative integer eigenvalues.
Furthermore, $I_{\infty}$ is semi-simple on $W$ with only
nonnegative integer eigenvalues and $H$ is semi-simple. \el

\begin{proof} For $p,q\in \Z$ with $p\ne q$, let $\g_{p,q}$ denote the
linear span of $E_{p,q}, E_{q,p}, E_{p,p}-E_{q,q}$, which is a
subalgebra isomorphic to $\sl_{2}$ with $E_{p,q}, E_{q,p},
E_{p,p}-E_{q,q}$ corresponding to $e,f,h$, respectively. With $W$ an
integrable $\gl_{\infty}$-module, we see that $W$ is an integrable
$\g_{p,q}$-module. Therefore, $E_{p,p}-E_{q,q}$ is semi-simple on
$W$ (by Remark \ref{rfd-integrable}).

Let $k$ be any fixed integer. We now show that $E_{k,k}$ is
semi-simple on $W$. Let $w$ be any vector of $W$. Then there exists
a finite subset $S$ of $\Z$ such that $k\in S$ and $E_{p,q}w=0$ for
all $p,q\in \Z$ with $q\notin S$. Consider the $\gl_{S}$-submodule
$U(\gl_{S})w$ generated by $w$. Suppose $u\in W$ satisfies that
$E_{p,q}u=0$ for all $p,q\in \Z$ with $q\notin S$. Then for any
$m,n\in S$ and $p,q\in \Z$ with $q\notin S$, we have
$$E_{p,q}E_{m,n}u=E_{m,n}E_{p,q}u+\delta_{q,m}E_{p,n}u-\delta_{n,p}E_{m,q}u
=E_{m,n}E_{p,q}u-\delta_{n,p}E_{m,q}u=0.$$
It follows from this and induction that $E_{p,q}U(\gl_{S})w=0$ for all $p,q\in \Z$ with $q\notin S$.
In particular, $E_{q,q}U(\gl_{S})w=0$ for all $q\notin S$.
Let $n\in S$. Pick an integer $q$ with $q\notin S$. We have $E_{q,q}=0$ on $U(\gl_{S})w$, so that
$E_{n,n}-E_{q,q}\ \left(=E_{n,n}\right)$ preserves $U(\gl_{S})w$.
Since $E_{n,n}-E_{q,q}$ is semi-simple on $W$ from the first paragraph,
$E_{n,n}-E_{q,q}$ is semi-simple on $U(\gl_{S})w$,  which implies that
$E_{n,n}$ is semi-simple on $U(\gl_{S})w$.
We have already proved that $E_{m,m}=0$ on $U(\gl_{S})w$ for $m\notin S$.
Thus $H$ preserves $U(\gl_{S})w$ and is semi-simple.
It then follows that $H$ is semi-simple on $W$.
In particular, $E_{k,k}$ is semi-simple on $W$.

Let $n\in \Z$ and let $u\in W$ be an eigenvector of $E_{n,n}$ with eigenvalue $\lambda\in \C$.
There exists an integer $m$ such that $m\ne n$ and $E_{p,m}u=0$ for all $p\in \Z$.  Then
$$(E_{n,n}-E_{m,m})u=E_{n,n}u=\lambda u \ \mbox{ and }
\ E_{n,m}u=0.$$ As $W$ is an integrable $\g_{m,n}$-module, in view
of Remark \ref{rsl2} we have $\lambda\in \N$ and
$(E_{m,n})^{\lambda+1}u=0$, as desired.
\end{proof}

Let $W$ be a $\gl_{\infty}$-module in category ${\cal{C}}_{int}$.
{}From Lemma \ref{lenn-eigenvalues} we have
$$W=\oplus _{\ell \in \N}W[\ell],$$ a direct sum of $\gl_{\infty}$-modules, where
$W[\ell]=\{ w\in W\;|\; I_{\infty}\cdot w=\ell w\}$ for $\ell\in
\N$. In view of this, it suffices to determine
$\gl_{\infty}$-modules in category ${\cal{C}}_{int}$, on which
$I_{\infty}$ acts as a nonnegative integer scalar.

As a refinement of Lemma \ref{lenn-eigenvalues} we have:

\bl{lsecond} Let $W$ be a $\gl_{\infty}$-module in category
${\cal{C}}_{int}$ such that $I_{\infty}$ acts as a nonnegative
integer scalar $\ell$. Then for every $n\in \Z$, $E_{n,n}$ is
semi-simple with nonnegative integer eigenvalues
 not exceeding $\ell$ and $(E_{m,n})^{\ell+1}=0$ for all $m,n\in \Z$ with $m\ne n$.
\el

\begin{proof} Let $m,n\in \Z$ with $m\ne n$.
Recall that $\g_{m,n}$ denotes the three-dimensional simple
subalgebra spanned by $E_{m,n},E_{n,m},E_{m,m}-E_{n,n}$. By Lemma
\ref{lenn-eigenvalues}, the eigenvalues of $E_{m,m}$ and $E_{n,n}$
are nonnegative integers which are bounded by $\ell$ as
$I_{\infty}=\ell$. It follows that the eigenvalues of
$E_{m,m}-E_{n,n}$ are bounded by $-\ell$ and $\ell$. As $W$ is an
integrable $\g_{m,n}$-module, it is a direct sum of irreducible
$\g_{m,n}$-modules of dimension bounded by $\ell+1$. Consequently,
$(E_{m,n})^{\ell+1}=0$.
\end{proof}

The following is an immediate consequence of Lemma \ref{lsecond}:

\bc{clevel-0} Let $W$ be a $\gl_{\infty}$-module in category
${\cal{C}}_{int}$ such that $I_{\infty}=0$. Then $\gl_{\infty}$ acts
trivially on $W$. \ec

Next, we have:

\bp{plevel-1} Let $W$ be a $\gl_{\infty}$-module in category
${\cal{C}}_{int}$ such that $I_{\infty}=1$ on $W$. Then $W$ is a
direct sum of submodules isomorphic to the natural module
$\C^{\infty}$. \ep

\begin{proof} Let $v$ be any nonzero $H$-weight vector. As $I_{\infty}=1$, there is $r\in \Z$ such that
$E_{r,r}v=v$ and $E_{n,n}v=0$ for all $n\in \Z$ with $n\ne r$. For
any $m\ne r$, we have
$$E_{m,m}E_{r,m}v=E_{r,m}E_{m,m}v+\delta_{m,r}E_{m,m}v-E_{r,m}v=-E_{r,m}v.$$
Since the only eigenvalues of $E_{r,r}$ are either $1$ or $0$ (by
Lemma \ref{lenn-eigenvalues}), we must have $E_{r,m}v=0$. Note that
$E_{m,m}v=0$ for $m\ne r$. Let $m,n\in \Z$ with $m\ne n$ and $m,n\ne
r$. We have $E_{r,r}E_{n,m}v=E_{n,m}v$ and
$E_{n,n}E_{n,m}v=E_{n,m}v$. As $I_{\infty}=1$, we must have
$E_{n,m}v=0$. Thus $E_{p,q}v=0$ for all $p,q\in \Z$ with $q\ne r$.

For $n\in \Z$, set $v_{n}=E_{n,r}v$. We have $v=v_{r}$. Assume $n\ne
r$. Then
\begin{eqnarray}
E_{n,n}v_{n}=E_{n,n}E_{n,r}v=E_{n,r}E_{n,n}v+E_{n,r}v=E_{n,r}v=v_{n}.
\end{eqnarray}
This forces $E_{m,m}v_{n}=0$ for all $m\ne n$ (no matter whether
$v_{n}\ne 0$). Next, we show that
\begin{eqnarray}\label{epqn=0}
E_{p,q}v_{n}=0\ \ \mbox{ for }p,q\in \Z\ \mbox{  with }\ q\ne n.
\end{eqnarray}
If $p=q\ne n$, we already have $E_{p,q}v_{n}=0$.
Assume $p\ne q,\ p\ne n,\ q\ne n$. We have
$$E_{q,q}(E_{p,q}v_{n})=E_{p,q}E_{q,q}v_{n}-E_{p,q}v_{n}
=-E_{p,q}v_{n}.$$
Again, since the only possible eigenvalues of $E_{q,q}$ are $0$ and $1$, we have $E_{p,q}v_{n}=0$.
Furthermore, we have $E_{m,r}v_{r}=E_{m,r}v=v_{m}$ and
\begin{eqnarray}\label{emnn=m}
E_{m,n}v_{n}=E_{m,n}E_{n,r}v=E_{n,r}E_{m,n}v+E_{m,r}v-\delta_{m,r}E_{n,n}v=E_{m,r}v=v_{m}
\end{eqnarray}
for $n\ne r$. It follows from (\ref{epqn=0}) and (\ref{emnn=m}) that
$U(\gl_{\infty})v$ is a homomorphism image of $\C^{\infty}$.
Consequently, $U(\gl_{\infty})v$ is isomorphic to $\C^{\infty}$.
Then it follows that $W$ is a sum of submodules isomorphic to
$\C^{\infty}$. As $\C^{\infty}$ is irreducible, $W$ is a direct sum
of submodules isomorphic to $\C^{\infty}$.
\end{proof}

Next, we shall prove that every $\gl_{\infty}$-module in category
${\cal{C}}_{int}$ is completely reducible. First, we prove a
technical result.

\bl{lpreparation} Let $W$ be a $\gl_{\infty}$-module in category
${\cal{C}}_{int}$. Suppose that $v$ is a nonzero $H$-eigenvector of
weight $\lambda$, satisfying the condition that
\begin{eqnarray}\label{eassumption-2.28}
&&E_{p,q}v=0\ \ \ \mbox{ for all }p,q\in \Z\ \mbox{ with }\ q\notin
\{r,r+1,\dots,n\},\nonumber\\
&&E_{r,j+1}v=0\ \ \ \mbox{ for }r\le j\le n,
\end{eqnarray}
where $r$ and $n$ are some fixed integers with $r<n$. Set
$$v'=(E_{r-1,r})^{\lambda_{r}-\lambda_{r+1}}\cdots
(E_{r-1,n-1})^{\lambda_{n-1}-\lambda_{n}}(E_{r-1,n})^{\lambda_{n}}v\in W.$$
Then
\begin{eqnarray}\label{eresult-2.29}
&&E_{p,q}v'=0\ \ \ \mbox{ for all }p,q\in \Z\ \mbox{ with }\ q\notin
\{r-1,r,\dots,n\},\nonumber\\
&&E_{r-1,j+1}v'=0\ \ \ \mbox{ for }r-1\le j\le n.
\end{eqnarray}
Furthermore, we have
\begin{eqnarray}
(E_{n,r-1})^{\lambda_{n}}(E_{n-1,r-1})^{\lambda_{n-1}-\lambda_{n}}\cdots
(E_{r,r-1})^{\lambda_{r}-\lambda_{r+1}}v'=\alpha v,
\end{eqnarray}
where
$\alpha=(\lambda_{n}!)^{2}((\lambda_{n-1}-\lambda_{n})!)^{2}\cdots
((\lambda_{r}-\lambda_{r+1})!)^{2}$, a nonzero integer.
 \el

\begin{proof} First of all, by Lemma \ref{lenn-eigenvalues}
we have $\lambda_{m}\in \N$ for all $m\in \Z$. {}From assumption
(\ref{eassumption-2.28}), $v$ is a singular vector in $W$ viewed as
a $\gl_{S}$-module with $S=\{ r,r+1,\dots,n\}$. Since $W$ is an
integrable $\gl_{S}$-module from assumption, we must have
$$\lambda_{j}-\lambda_{j+1}\in \N\ \ \ \mbox{ for }r\le j\le n-1.$$
Note that $\lambda_{m}=0$ for $m\notin S$ from
(\ref{eassumption-2.28}). In particular,
$\lambda_{r-1}=0=\lambda_{n+1}$.

Let $p,q\in \Z$ with $q\notin \{r-1,r,\dots,n\}$. For $r\le j\le n$,
we have
$$E_{p,q}E_{r-1,j}=E_{r-1,j}E_{p,q}-\delta_{p,j}E_{r-1,q}.$$
Then it follows from induction that
$$E_{p,q}v'=0\ \ \ \mbox{ for all }p,q\in \Z\ \mbox{with } q\notin
\{r-1,r,\dots,n\},$$ proving the first part of (\ref{eresult-2.29}).

For $r\le m\le n$, set
$$v_{m}=(E_{r-1,m})^{\lambda_{m}-\lambda_{m+1}}\cdots
(E_{r-1,n-1})^{\lambda_{n-1}-\lambda_{n}}(E_{r-1,n})^{\lambda_{n}}v.$$
Note that $v'=v_{r}$. Next, in several steps we prove that for every
$m>r$,
$$E_{r-1,j}v_{m}=0\ \ \ \mbox{ for }j=m,m+1,\dots,n,$$
{}from which we obtain the second part of (\ref{eresult-2.29}) by
taking $m=r$.

First, we have
\begin{eqnarray}\label{efirst-condition}
E_{j,m}(E_{r-1,m_{1}}\cdots E_{r-1,m_{s}}v)=0
\end{eqnarray}
for any $j,m,m_{1},\dots,m_{s}\in \Z$ with $r\le j<m,
m_{1},\dots,m_{s}$. This is because $E_{j,m}$ commutes with
$E_{r-1,m_{i}}$ for $i=1,\dots,s$ and $E_{j,m}v=0$ from assumption
(\ref{eassumption-2.28}).

Second, we have
\begin{eqnarray}\label{econ2.33}
E_{j,r-1}(E_{r-1,m_{1}}\cdots E_{r-1,m_{s}}v)=0
\end{eqnarray}
for any $j,m_{1},\dots,m_{s}\in \Z$ with $r\le j<m_{1},\dots,m_{s}$.
This follows from induction and (\ref{efirst-condition}), as for
$1\le i\le s$,
$$E_{j,r-1}E_{r-1,m_{i}}=E_{r-1,m_{i}}E_{j,r-1}+E_{j,m_{i}}$$ and
$E_{j,r-1}v=0$ by assumption (\ref{eassumption-2.28}).

Third, for any $j,m_{1},\dots,m_{s}\in \Z$ with $r\le
j<m_{1},\dots,m_{s}$, we have
\begin{eqnarray}\label{econ2.34}
(E_{j,j}-E_{r-1,r-1})(E_{r-1,m_{1}}\cdots
E_{r-1,m_{s}}v)=(\lambda_{j}-s)(E_{r-1,m_{1}}\cdots E_{r-1,m_{s}}v),
\end{eqnarray}
noticing that $E_{j,j}v=\lambda_{j}v$ and $E_{r-1,r-1}v=0$ by
assumption (\ref{eassumption-2.28}).

Fourth, for any $j,m_{1},\dots,m_{s}\in \Z$ with $r\le
j<m_{1},\dots,m_{s}$, we have
\begin{eqnarray}\label{e401}
\left(E_{r-1,j}\right)^{(\lambda_{j}-s)+1}(E_{r-1,m_{1}}\cdots
E_{r-1,m_{s}}v)=0
\end{eqnarray}
and
\begin{eqnarray}\label{e402}
&&\left(E_{j,r-1}\right)^{\lambda_{j}-s}\left(E_{r-1,j}\right)^{\lambda_{j}-s}(E_{r-1,m_{1}}\cdots
E_{r-1,m_{s}}v)\nonumber\\
&=&((\lambda_{j}-s)!)^{2}(E_{r-1,m_{1}}\cdots E_{r-1,m_{s}}v).
\end{eqnarray}
Recall that for a fixed $j$ with the above condition, $E_{j,r-1}$,
$E_{r-1,j}$, $E_{j,j}-E_{r-1,r-1}$ span a Lie subalgebra isomorphic
to $\sl(2,\C)$. Then the two assertions follow from (\ref{econ2.33})
and (\ref{econ2.34}) as $W$ is an integrable $\sl(2,\C)$-module from
assumption.

Fifth, we apply this to $v_{m}$ for $r\le m\le n$. In view of
(\ref{e401}) we have
\begin{eqnarray}\label{e501}
E_{r-1,m}v_{m}=(E_{r-1,m})^{(\lambda_{m}-\lambda_{m+1})+1}v_{m+1}=0,
\end{eqnarray}
where we set $v_{n+1}=v$. As $E_{r-1,j}$ with $j\ge r$ commute each
other, by (\ref{e501}) we obtain
\begin{eqnarray}
E_{r-1,j}v_{r}=0\ \ \ \mbox{ for }r\le j\le n.
\end{eqnarray}
This also holds  for $j>n$ since $E_{r-1,j}v=0$ from assumption.
Now we have proved the second part of (\ref{eresult-2.29}). The last
assertion follows immediately from repeatedly applying (\ref{e402}).
\end{proof}

On the basis of Lemma \ref{lpreparation} we have:

\bp{psingular-vector} Let $W$ be an integrable $\gl_{\infty}$-module
and let $v\in W$ be a nonzero $H$-weight vector satisfying the
condition that
\begin{eqnarray}
&&E_{p,q}v=0\ \ \ \mbox{ for }p,q\in \Z \ \mbox{ with }\ q\notin \{r,r+1,\dots,n\},\nonumber\\
&&E_{j,j+1}v=0\ \ \ \mbox{ for }r\le j\le n-1,
\end{eqnarray}
where $r$ and $n$ are some integers with $r<n$. Then the submodule
$U(\gl_{\infty})v$ is irreducible. \ep

\begin{proof} Let $u$ be any nonzero vector in $U(\gl_{\infty})v$. We now prove
$v\in U(\gl_{\infty})u$, so that
$U(\gl_{\infty})u=U(\gl_{\infty})v$. As $u\in U(\gl_{\infty})v$,
there exist two integers $r'$ and $n'$ with $r'\le r<n\le n'$ such
that $u\in U(\gl_{S'})v$, where $$S'=\{ r',r'+1,\dots,n'\}.$$ Note
that if $U(\gl_{S'})v$ is an irreducible $\gl_{S'}$-module, then we
have $v\in U(\gl_{S'})v=U(\gl_{S'})u\subset U(\gl_{\infty})u$, as
desired.

As $W$ is an integrable $\gl_{S'}$-module, any $\gl_{\S'}$-submodule
of $W$ generated by a singular vector is irreducible. In view of
this, it suffices to prove that there is a singular vector $v'$ such
that $U(\gl_{S'})v=U(\gl_{S'})v'$. Note that from assumption we have
$E_{r,j}v=0$ for $r+1\le j\le n$ and $E_{r,q}v=0$ for $q\ge n+1$.
Thus $E_{r,j}v=0$ for all $j\ge r+1$. To summarize we have
\begin{eqnarray*}
&&E_{p,q}v=0\ \ \ \mbox{ for all }p,q\in \Z\ \mbox{ with } q\notin
\{r,r+1,\dots,n'\},\\
 &&E_{r,j}v=0\ \ \ \mbox{ for }j\ge r+1.
\end{eqnarray*}
By repeatedly applying Lemma \ref{lpreparation}, we obtain a vector
$v'$ such that
\begin{eqnarray*}
&&E_{p,q}v=0\ \ \ \mbox{ for all }p,q\in \Z\ \mbox{ with } q\notin
\{r',r'+1,\dots,n'\},\\
 &&E_{r',j}v=0\ \ \ \mbox{ for }j\ge r'+1
\end{eqnarray*}
and such that $U(\gl_{S'})v=U(\gl_{S'})v'$, where the second
condition implies that $v'$ is a singular vector of the
$\gl_{S'}$-module $U(\gl_{S'})v$. This proves that $U(\gl_{S'})v$ is
an irreducible $\gl_{S'}$-module,
$U(\gl_{\infty})u=U(\gl_{\infty})v$ for any nonzero vector $u$ in
$U(\gl_{\infty})v$, and $U(\gl_{\infty})v$ is an irreducible
$\gl_{\infty}$-module.
\end{proof}

Now, we are in a position to present our main result.

\bt{tmain-D} Every $\gl_{\infty}$-module in category
${\cal{C}}_{int}$ is completely reducible.  \et

\begin{proof}
Let $W$ be any $\gl_{\infty}$-module in category ${\cal{C}}_{int}$.
By Lemma \ref{lenn-eigenvalues}, $H$ is semi-simple on $W$. To prove
$W$ is completely reducible, it suffices to show that the submodule
generated by each vector is a sum of irreducible submodules.

Let $w\in W$. We shall prove that $U(\gl_{\infty})w$ is a sum of
irreducible submodules. {}From assumption, there are integers $r$
and $n$ with $r<n$ such that
$$E_{p,q}w=0\ \ \ \mbox{ for all }p,q\in \Z\ \mbox{with }q\notin
\{r,r+1,\dots,n\}.$$ Set $S=\{r,r+1,\dots,n\}.$ As a
$\gl_{S}$-submodule of $W$, $U(\gl_{S})w$ is integrable, so that
$U(\gl_{S})w$ is a direct sum of finite-dimensional irreducible
$\gl_{S}$-modules. Now, it suffices to show that for every
finite-dimensional irreducible $\gl_{S}$-submodule $U$ of
$U(\gl_{S})w$, $U(\gl_{\infty})U$ is irreducible. Let $v\in U$ be a
highest weight vector of weight $\lambda$. Then
$$U(\gl_{\infty})U=U(\gl_{\infty})v.$$
Now we prove that $U(\gl_{\infty})v$ is irreducible. Note that for
any $p,q\in \Z$ with $q\notin S$ and for any $m,k\in S$,
$$E_{p,q}E_{m,k}=E_{m,k}E_{p,q}-\delta_{p,k}E_{m,q}.$$
By induction we get
$$E_{p,q}U(\gl_{S})w=0\ \ \ \mbox{ for all }p,q\in \Z\ \mbox{with }q\notin
S.$$
In particular, we have
$$E_{p,q}v=0\ \ \ \mbox{ for all }p,q\in \Z\ \mbox{with }q\notin
S=\{r,r+1,\dots,n\}.$$ As $v$ is a highest weight vector, we also
have
$$E_{j,j+1}v=0\ \  \mbox{ for }\ r\le j\le n-1.$$
By Proposition \ref{psingular-vector}, $U(\gl_{\infty})v$ is an
irreducible $\gl_{\infty}$-module.
 Therefore, $W$ is completely reducible.
\end{proof}

\section{Classification of irreducible $\gl_{\infty}$-modules in ${\cal{C}}$ and ${\cal{C}}_{int}$}

In this section, we classify irreducible $\gl_{\infty}$-modules in
categories ${\cal{C}}$ and ${\cal{C}}_{int}$. To achieve this goal,
for any finite subset $S$ of $\Z$ and for any irreducible
$\gl_{S}$-module $U$, through a generalized Verma module
construction we construct an irreducible $\gl_{\infty}$-module
$L(S,U)$ in category ${\cal{C}}$ and we show that any irreducible
module in ${\cal{C}}$ is isomorphic to a module of this form.
Furthermore, for a linear functional $\lambda$ on $H$ compatible
with $S$ in a certain sense we define a generalized Verma
$\gl_{\infty}$-module $M(S,\lambda)$ and construct an irreducible
module $L(S,\lambda)$. We then determine when $L(S,\lambda)$ is
integrable and we show that every irreducible module in
${\cal{C}}_{int}$ is isomorphic to such an integrable module
$L(S,\lambda)$.

Let $S$ be a {\em finite} and {\em nonempty} subset of $\Z$, which
is fixed temporarily. Recall
\begin{eqnarray*}
\gl_{S}={\rm span}\{ E_{m,n}\ |\ m,n\in S\}\subset \gl_{\infty}.
\end{eqnarray*}
To $S$, we associate a triangular decomposition
\begin{eqnarray}
\gl_{\infty}=\gl_{\infty}^{(S,+)}\oplus \gl_{\infty}^{(S,0)}\oplus
\gl_{\infty}^{(S,-)},
\end{eqnarray}
where
$\gl_{\infty}^{(S,0)}=\gl_{S}$,
\begin{eqnarray}
 \gl_{\infty}^{(S,+)}&=&\span\{ E_{m,n}\;|\; m,n\in
\Z,\; n\notin S\},\nonumber\\
\gl_{\infty}^{(S,-)}&=&\span\{ E_{m,n}\;|\; m,n\in \Z,\; m\notin
S,\;  n\in S\}.
\end{eqnarray}
Notice that $\gl_{\infty}^{(S,-)}$ is an abelian subalgebra and that
$\gl_{\infty}^{(S,+)}+\gl_{S}$ is a semi-product.

Recall that
\begin{eqnarray*}
I_{S}=\sum_{n\in S}E_{n,n}\in \gl_{S}.
\end{eqnarray*}
We have
\begin{eqnarray}
&&[I_{S},E_{m,n}]=0\ \mbox{ for either }m,n\in S, \mbox{ or }m,n\notin S,\nonumber\\
&&[I_{S},E_{m,n}]=E_{m,n}\ \mbox{ for }m\in S,\ n\notin S,\nonumber\\
&&[I_{S},E_{m,n}]=-E_{m,n}\ \mbox{ for }m\notin S,\; n\in S.
\end{eqnarray}
Using $\ad(I_{S})$, we make $\gl_{\infty}$ a $\Z$-graded Lie algebra
for which
\begin{eqnarray}
&&\gl_{\infty}^{(n)}=0\ \ \mbox{ if }|n|\ge 2,\nonumber\\
&&\gl_{\infty}^{(0)}={\rm span} \{E_{m,n}\ |\ \mbox{ either }m,n\in S, \mbox{ or }m,n\notin S\},\nonumber\\
&&\gl_{\infty}^{(1)}={\rm span} \{E_{m,n}\ |\ m\in S,\ n\notin S\},\nonumber\\
&&\gl_{\infty}^{(-1)}={\rm span} \{E_{m,n}\ |\ m\notin S,\ n\in S\}.
\end{eqnarray}

Set
\begin{eqnarray}
\gl_{S^{o}}={\rm span}\{ E_{m,n}\ |\ m,n\notin S\}.
\end{eqnarray}
Then
\begin{eqnarray}
\gl_{\infty}^{(0)}=\gl_{S}\oplus \gl_{S^{o}},
\end{eqnarray}
a direct product. Notice that $\gl_{\infty}^{(\pm 1)}$ both are
abelian subalgebras. We see that
\begin{eqnarray}
\gl_{\infty}^{(S,+)}=\gl_{\infty}^{(1)}\oplus \gl_{S^{o}},\ \
\ \ \gl_{\infty}^{(S,-)}=\gl_{\infty}^{(-1)}.
\end{eqnarray}

Let $U$ be a $\gl_{S}$-module. Let $\gl_{\infty}^{(S,+)}$ act
trivially on $U$, to make $U$ a
$(\gl_{\infty}^{(S,+)}+\gl_{S})$-module. Then form a generalized
Verma $\gl_{\infty}$-module
\begin{eqnarray}
M(S,U)=U(\gl_{\infty})\otimes
_{U\left(\gl_{\infty}^{(S,+)}+\gl_{S}\right)}U.
\end{eqnarray}
In view of the P-B-W theorem we have
$$M(S,U)=U(\gl_{\infty}^{(-1)})\otimes U=S(\gl_{\infty}^{(-1)})\otimes
U.$$ By endowing $U$ with degree $0$, we make $M(S,U)$ a $\Z$-graded
$\gl_{\infty}$-module. (The homogeneous subspaces are
infinite-dimensional in general.)

\bl{lM(U)} The $\gl_{\infty}$-module $M(S,U)$ belongs to the
category ${\cal{C}}$. \el

\begin{proof} Let $W$ consist of every $w\in M(S,U)$, satisfying the
condition that there exists a finite subset $T$ of $\Z$ such that
$E_{m,n}w=0$ for all $m,n\in \Z$ with $n\notin T$. By the
construction we have $U\subset W$. Then it suffices to prove that
$W$ is a submodule. Assume that $w\in W$ with a finite subset $T$ of
$\Z$ such that $E_{m,n}w=0$ for all $m,n\in \Z$ with $n\notin T$.
Let $p,q\in \Z$ be arbitrarily fixed. For $m,n\in \Z$, we have
$$E_{m,n}(E_{p,q}w)=E_{p,q}E_{m,n}w+\delta_{n,p}E_{m,q}w-\delta_{q,m}E_{p,n}w.$$
One sees that $E_{m,n}(E_{p,q}w)=0$ for all $m,n\in\Z$ with $n\notin
T\cup \{ p\}$. This proves $E_{p,q}w\in W$. Then the lemma follows.
\end{proof}

\bd{dL(S,U)} {\em Let $U$ be a $\gl_{S}$-module as before. Denote by $L(S,U)$ the
quotient of the $\Z$-graded $\gl_{\infty}$-module $M(S,U)$ by the
maximal graded submodule with trivial degree-$0$ homogeneous
subspace.}
\ed

\br{rVerma-grading} {\em Assume that $U$ is a $\gl_{S}$-module on
which $I_{S}$ acts as a scalar $\alpha\in \C$. Then $M(S,U)$ is a
canonically graded $\gl_{\infty}$-module
\begin{eqnarray}
M(S,U)=\bigoplus _{n\in \N}M(S,U)_{\alpha-n},
\end{eqnarray}
where $M(S,U)_{\alpha-n}=\{ w\in M(S,U)\; |\; I_{S}\cdot
w=(\alpha-n)w\}$ for $n\in \N$.} \er

Let $W$ be a $\gl_{\infty}$-module and let $S$ be a finite subset of
$\Z$ as before. Set
\begin{eqnarray}
\Omega_{S}(W)=\left\{ w\in W\; |\; E_{p,q}w=0\ \ \mbox{ for all
}p,q\in \Z\ \mbox{ with }q\notin S\right\}.
\end{eqnarray}
It can be readily seen that $\Omega_{S}(W)$ is a
$\gl_{S}$-submodule of $W$.

\bp{panyC} Let $S$ be a finite subset of $\Z$ and let $U$ be an
irreducible $\gl_{S}$-module. Then $L(S,U)$ is an irreducible $\gl_{\infty}$-module
belong to category ${\cal{C}}$. On the other hand, every irreducible
$\gl_{\infty}$-module in category ${\cal{C}}$ is isomorphic to
$L(S,U)$ for some finite subset $S$ of $\Z$ and for
some irreducible $\gl_{S}$-module $U$. \ep

\begin{proof} Since $U$ is an
irreducible $\gl_{S}$-module, $U$ is necessarily
countable-dimensional. Then $I_{S}$ acts on $U$ as a scalar, say
$\alpha\in \C$. From Remark \ref{rVerma-grading}, $M(S,U)$ is
canonically $\C$-graded by the eigenspaces of $I_{S}$, where $M(S,U)=\oplus_{n\in
\N}M(S,U)_{\alpha-n}$ with $M(S,U)_{\alpha}=U$. We see that every
submodule of $M(S,U)$ is graded. It follows that $M(S,U)$ has a
unique maximal submodule. Consequently, $L(S,U)$ is an irreducible
$\gl_{\infty}$-module. From Lemma \ref{lM(U)}, $L(S,U)$ belongs to ${\cal{C}}$.

Now, let $W$ be an irreducible $\gl_{\infty}$-module in category
${\cal{C}}$. Pick a nonzero vector $w$ in $W$. Then there exists a
finite and nonempty subset $S$ of $\Z$ such that $w\in
\Omega_{S}(W)$.  Set $U=U(\gl_{S})w$, a $\gl_{S}$-submodule of
$\Omega_{S}(W)$. Using the P-B-W theorem we get
$$W=U(\gl_{\infty})w=U(\gl_{\infty}^{(-1)})U.$$
As $W$ belongs to category ${\cal{C}}$, $I_{\infty}$ acts on $W$ and
commutes with the action of $\gl_{\infty}$. With $W$ an irreducible
$\gl_{\infty}$-module, $W$ must be countable dimensional (over
$\C$). It then follows that $I_{\infty}$ acts as a scalar on $W$,
say $\alpha$. Notice that $I_{S}\cdot w=I_{\infty}\cdot w=\alpha w$.
Consequently, $I_{S}$ acts on $U$ as scalar $\alpha$. As
$W=U(\gl_{\infty}^{(-1)})U$, $I_{S}$ is semisimple on $W$ with
eigenvalues contained in $\{\alpha-n\ |\ n\in \N\}$ and $U$ is the
eigenspace of eigenvalue $\alpha$. It follows that $U$ is an
irreducible $\gl_{S}$-module. By the construction of $M(S,U)$, there exists an epimorphism from
$M(S,U)$ to $W$, which reduces to an isomorphism from
$L(S,U)$ to $W$. This completes the proof.
\end{proof}

From the second part of the proof of Proposition \ref{panyC} we immediately have:

\bl{lSvacuum}
Let $W$ be an irreducible $\gl_{\infty}$-module in ${\mathcal{C}}$. Then there exists
a finite subset $S$ of $\Z$ such that $\Omega_{S}(W)\ne 0$. Furthermore, for any
such finite subset $S$ of $\Z$,
$\Omega_{S}(W)$ is an irreducible $\gl_{S}$-module and $W\simeq L(S,\Omega_{S}(W))$.
\el

The following is also immediate:

\bl{lvacuum-irred}
Let $S$ be a finite subset of $\Z$ and let $U$ be an irreducible $\gl_{S}$-module.
Then $\Omega_{S}(L(S,U))=U$.
\el

We next determine the isomorphism classes of irreducible
$\gl_{\infty}$-modules $L(S,U)$. Let $S$ be a finite subset of $\Z$
and $S_{1}$ a subset of $S$. Assume that $U_{1}$ is an irreducible
$\gl_{S_{1}}$-module on which $I_{S_{1}}$ acts as a scalar
$\alpha\in \C$. In the following we associate an irreducible
$\gl_{S}$-module to $U_{1}$. Set
$$N={\rm span}\{ E_{p,q}\ |\ p,q\in S,\; q\notin S_{1}\}\ \ \mbox{ and }\ \
B=N+\gl_{S_{1}}.$$ We see that both $B$ and $N$ are subalgebras of
$\gl_{S}$ and $B$ contains $N$ as an ideal. Letting $N$ act on
$U_{1}$ trivially, we make $U_{1}$ a $B$-module. Then form an
induced module
\begin{eqnarray}
{\rm
ind}_{\gl_{S_{1}}}^{\gl_{S}}(U_{1})=U(\gl_{S})\otimes_{U(B)}U_{1}.
\end{eqnarray}
As we have seen before, $I_{S_{1}}$ gives rise to an
$(\alpha+\Z)$-grading on ${\rm ind}_{\gl_{S_{1}}}^{\gl_{S}}(U_{1})$
with $U_{1}$ as the degree-$\alpha$ subspace. It follows that ${\rm
ind}_{\gl_{S_{1}}}^{\gl_{S}}(U_{1})$ has a unique maximal submodule.
Then we define $U_{1}^{S}$ to be the (unique) irreducible quotient
$\gl_{S}$-module of ${\rm ind}_{\gl_{S_{1}}}^{\gl_{S}}(U_{1})$.

\bl{lu1us} Let $S$ be a finite subset of $\Z$ and $S_{1}$ a subset
of $S$. Assume that $U_{1}$ is an irreducible $\gl_{S_{1}}$-module.
Then $L(S_{1},U_{1})\simeq L(S,U_{1}^{S})$. \el

\begin{proof} Set $U=U(\gl_{S})U_{1} \subset L(S_{1},U_{1}).$ Note that
$E_{p,q}U_{1}=0$ for $p,q\in \Z$ with $q\notin S$. It follows that
$$E_{p,q}\cdot U(\gl_{S})U_{1}=0\ \ \ \mbox{ for }
p,q\in \Z\ \mbox{with }q\notin S.$$ Then
$$L(S_{1},U_{1})=U(\gl_{\infty})\cdot U$$
and $I_{S}$ acts on $U=U(\gl_{S})U_{1}$ also as scalar $\alpha$
since $I_{S}=I_{S_{1}}$ on $U_{1}$. Just as with $M(S,U)$, we see
that $L(S_{1},U_{1})$ is naturally an $(\alpha+\Z)$-graded
$\gl_{\infty}$-module by $I_{S}$, with $U$ as the degree-$\alpha$
subspace. Since $L(S_{1},U_{1})$ is irreducible, it follows that $U$
is an irreducible $\gl_{S}$-module and $L(S_{1},U_{1})\simeq
L(S,U)$. {}From the construction of $U_{1}^{S}$, we have $U\simeq
U_{1}^{S}$. Thus $L(S_{1},U_{1})\simeq L(S,U_{1}^{S})$.
\end{proof}

As an immediate consequence of Lemmas \ref{lvacuum-irred} and \ref{lu1us} we have:

\bc{cu1u2} Let $S_{1}$ and $S_{2}$ be finite subsets of $\Z$ and let
$U_{1}$ and $U_{2}$ be irreducible modules for $\gl_{S_{1}}$ and
$\gl_{S_{2}}$, respectively. Set $S=S_{1}\cup S_{2}$. Then
$L(S_{1},U_{1})\simeq L(S_{2},U_{2})$ if and only if
$U_{1}^{S}\simeq U_{2}^{S}$. \ec

Next, we study an analog of Verma module. For $i\in \Z$, let
$\varepsilon_{i}$ denote the linear functional on $H$ defined by
$\varepsilon_{i}(E_{j,j})=\delta_{i,j}$ for $j\in \Z$. The root
system of $\gl_{\infty}$ with respect to Cartan subalgebra $H$ is
given by
\begin{eqnarray}
\Delta=\{ \varepsilon_{i}-\varepsilon_{j}\ |\ i,j\in \Z,\ i\ne j\}.
\end{eqnarray}
Note that the usual polarization is given by $\Delta_{\pm}=\{\pm
(\varepsilon_{i}-\varepsilon_{j})\ |\ i,j\in \Z,\ i<j\}$.

Let $S$ be a subset of $\Z$. Recall $\gl_{S}={\rm span}\{E_{m,n}\ |\
m,n\in S\}$. Set
\begin{eqnarray}
\gl_{S}^{\pm}&=&{\rm span}\{ E_{m,n}\ |\ m,n\in S,\ \pm
(n-m)>0\},\nonumber\\
\gl_{S}^{0}&=&H_{S}={\rm span}\{ E_{n,n}\ |\ n\in S\}.
\end{eqnarray}
 For $\lambda\in H^{*}$, set
$${\rm supp}(\lambda)=\{ m\in \Z\ |\ \lambda_{m}\ (=\lambda(E_{m,m}))\ne 0\}.$$
Let $\lambda\in H^{*}$ with ${\rm supp}(\lambda)\subset S$. Denote
by $\lambda_{S}$ the restriction of $\lambda$ on $H_{S}$. Let
$M(\lambda_{S})$ and $L(\lambda_{S})$ denote the Verma module and
the irreducible quotient module for Lie algebra $\gl_{S}$,
respectively.

Define $M(S,\lambda)$ to be the generalized Verma
$\gl_{\infty}$-module $M(S,U)$ with $U=M(\lambda_{S})$. On the other
hand, we have a generalized Verma $\gl_{\infty}$-module
$M(S,L(\lambda_{S}))$, which is a quotient module of $M(S,\lambda)$.
Furthermore, denote by $L(S,\lambda)$ the irreducible quotient
module of $M(S,L(\lambda_{S}))$.

Set
\begin{eqnarray}
\Delta_{-}(S)=\{(\varepsilon_{i}-\varepsilon_{j})\ |\ i,j\in S,\
i>j\}\cup \{ (\varepsilon_{p}-\varepsilon_{i})\ |\ i\in S,\ p\notin
S\}.
\end{eqnarray}
Furthermore, set
\begin{eqnarray}
Q_{-}(S)=\N\cdot \Delta_{-}(S)\subset H^{*}.
\end{eqnarray}
We have
\begin{eqnarray}
M(S,\lambda)=\oplus_{\alpha\in
Q_{-}(S)}M(S,\lambda)_{\lambda+\alpha},
\end{eqnarray}
where $M(S,\lambda)_{\lambda}$ is $1$-dimensional. It is
straightforward to show that every weight subspace is
finite-dimensional.

\bd{dsingular-vector} {\em Let $W$ be a $\gl_{\infty}$-module and
let $S$ be a subset of $\Z$. A nonzero vector $v\in W$ is called an
{\em $S$-singular vector} if $v$ is an $H$-eigenvector such that
\begin{eqnarray}
&&E_{p,q}v=0\ \ \ \mbox{ for }p,q\in \Z \ \mbox{ with }\ q\notin S,\nonumber\\
&&E_{m,n}v=0\ \ \ \mbox{ for }m,n\in S\ \mbox{ with }m<n.
\end{eqnarray}}
\ed

The following universal property of $M(S,\lambda)$ is
straightforward to prove:

\bl{luniversal} Let $W$ be a $\gl_{\infty}$-module and let $w$ be an
$S$-singular vector of weight $\lambda$ in $W$. Then there exists a
$\gl_{\infty}$-module homomorphism $\theta$ from $M(S,\lambda)$ to
$W$, uniquely determined by $\theta(v)=w$, where $v$ is an
$S$-singular vector in $M(S,\lambda)$ of weight $\lambda$.
Furthermore, if $W$ is irreducible, we have $W\simeq
L(S,\lambda)$.\el

We also have the following result:

\bl{lsingular-vector} Let $S$ be a finite subset of $\Z$ and let
$\lambda\in H^{*}$ with ${\rm supp}(\lambda)\subset S$. Then
$S$-singular vectors in $L(S,\lambda)$ are unique up to scalar
multiples. \el

\begin{proof} Set $\ell=\sum_{m\in \Z}\lambda_{m}=\sum_{m\in S}\lambda_{m}.$
Then $I_{S}\cdot w=\ell w$ for $w\in L(\lambda_{S})$ (the
irreducible highest weight $\gl_{S}$-module of highest weight
$\lambda_{S}$). We see that $L(S,\lambda)$ is an $(\ell+\Z)$-graded
$\gl_{\infty}$-module with
$$L(S,\lambda)=\bigoplus_{n\in \N} L(S,\lambda)_{\ell-n},$$ where
$L(S,\lambda)_{\ell-n}=\{ w\in L(S,\lambda)\; |\; I_{S}\cdot
w=(\ell-n)w\}$ and $L(S,\lambda)_{\ell}=L(\lambda_{S})$. It is
straightforward to see that the subspace spanned by all $S$-singular
vectors in $L(S,\lambda)$ is $I_{S}$-stable, so that it is a graded
subspace. Let $u$ be any homogeneous $S$-singular vector of degree
$\ell-k$ with $k\in \N$. We have
$$U(\gl_{\infty})u=U(\gl_{\infty}^{(-1)})U(\gl_{S})u
\subset \bigoplus_{n\ge k}L(S,\lambda)_{\ell-n}.$$ Since
$L(S,\lambda)$ is irreducible, $k$ must be zero, so that every
$S$-singular vector is contained in $L(S,\lambda)_{\ell}$. Then each
$S$-singular vector is a singular vector in $L(\lambda_{S})$ viewed
as a $\gl_{S}$-module, which is known to be unique up to scalar
multiples. Consequently, $S$-singular vectors in $L(S,\lambda)$ are
unique up to scalar multiples.
\end{proof}

As an immediate consequence we have:

\bc{cclassification-0} Let $S$ be a finite subset of $\Z$ and let
$\lambda,\mu\in H^{*}$ such that ${\rm supp}(\lambda)$, ${\rm
supp}(\mu)\subset S$. Then $L(S,\lambda)\simeq L(S,\mu)$ if and only
if $\lambda=\mu$. \ec

\bl{lany-verma} Let $W$ be an irreducible $\gl_{\infty}$-module in
category ${\cal{C}}$, satisfying the condition that for any vector
$w\in W$ and for any finite subset $S$ of $\Z$, $U(\gl_{S})w$ is a
$\gl_{S}$-module in category $o$. Then $W\simeq L(S,\lambda)$ where
$S=\{ r,r+1,\dots,n\}$ for some integers $r$ and $n$ with $r<n$ and
for some $\lambda\in H^{*}$ with ${\rm supp}(\lambda)\subset S$. \el

\begin{proof} Let $w$ be a nonzero vector in $W$. As $W$ belongs to category ${\cal{C}}$,
 there exist integers $r$ and $n$ with $r<n$ such that
$$E_{p,q}w=0\ \ \ \mbox{ for all }p,q\in \Z\ \mbox{ with }q\notin \{r,r+1,\dots,n\}.$$
Set $S=\{r,r+1,\dots,n\}$. We have $w\in \Omega_{S}(W)$ and hence
$U(\gl_{S})w\subset \Omega_{S}(W)$. {}From our assumption, in the
$\gl_{S}$-submodule $U(\gl_{S})w$ there exists a highest weight
vector $v$. Then we have
\begin{eqnarray*}
&&E_{p,q}v=0\ \ \ \mbox{ for all }p,q\in \Z \mbox{ with } q\notin
S,\\
&&E_{m,n}v=0\ \ \ \mbox{ for }m,n\in S \mbox{ with } m<n.
\end{eqnarray*}
On the other hand, as $E_{q,q}v=0$ for $q\in \Z\backslash S$, $v$ is
an $H$-eigenvector of some weight $\lambda\in H^{*}$ with ${\rm
supp}(\lambda)\subset S$. Then $v$ is an $S$-singular vector of
weight $\lambda$. It follows from Lemma \ref{luniversal} that
$W\simeq L(S,\lambda)$.
\end{proof}

By Corollary \ref{cu1u2} we have:

\bc{cisomorphism} Let $S_{1},S_{2}$ be finite subsets of $\Z$ and
let $\lambda,\mu\in H^{*}$ be such that ${\rm supp}(\lambda)\subset
S_{1}$ and ${\rm supp}(\mu)\subset S_{2}$. Then
$L(S_{1},\lambda)\simeq L(S_{2},\mu)$ if and only if
$L(\lambda_{S_{1}})^{S}\simeq L(\mu_{S_{2}})^{S}$, where
$S=S_{1}\cup S_{2}$, and $L(\lambda_{S_{1}})$, $L(\mu_{S_{2}})$ are
the irreducible highest weight modules for $\gl_{S_{1}}$ and
$\gl_{S_{2}}$, respectively.\ec

\br{rautomorphism} {\em Let $\sigma$ be a permutation on $\Z$. It
can be readily seen that $\sigma$ becomes an automorphism of the Lie
algebra $\gl_{\infty}$ by defining
\begin{eqnarray}
\sigma(E_{m,n})=E_{\sigma(m),\sigma(n)}\ \ \ \mbox{ for }m,n\in \Z.
\end{eqnarray}
For any $\gl_{\infty}$-module $W$, we denote by $W^{[\sigma]}$ the
$\gl_{\infty}$-module with $W$ as the underlying space and with the
action given by
$$a\cdot w=\sigma(a)w\ \ \ \mbox{ for }a\in \gl_{\infty},\ w\in W.$$
We see that if $W$ is an integrable $\gl_{\infty}$-module, then
$W^{[\sigma]}$ is still an integrable module. Furthermore, if $W$ is
in category ${\cal{C}}_{int}$, $W^{[\sigma]}$ is still in category
${\cal{C}}_{int}$.} \er

The following is straightforward to prove:

\bl{lshift} Let $\lambda\in H^{*}$ be such that ${\rm
supp}(\lambda)\subset S$, and let $\sigma$ be a permutation on $\Z$
such that $\sigma(m)<\sigma(n)$ for any $m,n\in S$ with $m<n$. Then
$L(S,\lambda)^{[\sigma]}\simeq L(\sigma^{-1}(S),\lambda\circ
\sigma)$. \el

Next, we determine when $L(S,\lambda)$ is an integrable module.

\bd{dp+} {\em Denote by $P_{+}(S)$ the set of $\lambda\in H^{*}$
such that
\begin{eqnarray}
&&\lambda(E_{i,i})\in \N\ \ \ \mbox{ for }i\in \Z,\nonumber\\
&&\lambda(E_{i,i})=0\ \ \ \mbox{ whenever }i\notin S,\nonumber\\
&&\lambda(E_{i,i})\ge \lambda(E_{j,j})\ \ \ \mbox{ for }i,j\in S\
\mbox{ with }i<j.
\end{eqnarray}}
\ed

Note that if $\lambda \in P_{+}(S)$, then ${\rm supp}(\lambda)\subset
S$.

\br{rsimple-fact} {\em We here mention a fact which we need in the
proof of the next proposition. Let $\lambda$ be a (dominant
integral) weight for Lie algebra $\gl_{n+1}$ such that
$$\lambda_{j}\in
\N\ \ \mbox{ for }1\le j\le n+1, \mbox{ and }\lambda_{1}\ge
\lambda_{2}\ge \cdots \ge \lambda_{n+1}.$$ Then the irreducible
highest weight $\gl_{n+1}$-module $L(\lambda)$ is
finite-dimensional. Let $\mu$ be the lowest weight of $L(\lambda)$.
We claim that $\mu_{j}\in \N$ for $1\le j\le n+1$. It was known (cf.
\cite{hum}) that $\mu=\sigma(\lambda)$ where $\sigma$ is the longest
Weyl group element. Suppose $\gamma$ is any weight such that
$\gamma_{j}\in \N$ for $1\le j\le n+1$. For $1\le i\le n$, with the
reflection $r_{i}$, we have
$$r_{i}(\gamma)=\gamma-\<\gamma,\alpha_{i}^{\vee}\>\alpha_{i}
=\gamma-(\gamma_{i}-\gamma_{i+1})(\varepsilon_{i}-\varepsilon_{i+1}).$$
Then
 \begin{eqnarray*}
r_{i}(\gamma)(E_{j,j})=\begin{cases}\gamma_{j} & \mbox{ if }j\ne i, i+1,\\
\gamma_{i+1}&\mbox{ if }j=i,\\
\gamma_{i}&\mbox{ if }j=i+1,
\end{cases}
\end{eqnarray*}
which implies $r_{i}(\gamma)_{j}\in \N$ for $1\le j\le n+1$. Then
the claim follows from induction.} \er

We  have:

\bp{pintegrability}  Let $S$ be a finite subset of $\Z$ and let
$\lambda\in H^{*}$ be such that ${\rm supp}(\lambda)\subset S$. Then
$L(S,\lambda)$ is an integrable $\gl_{\infty}$-module if and only if
$\lambda\in P_{+}(S)$. \ep

\begin{proof} In view of Remark \ref{cisomorphism} and Lemma \ref{lshift},
it suffices to prove the proposition for $S=\{
r,r+1,\dots,n\}\subset \Z$, where $r$ and $n$ are fixed integers
with $r<n$.

Assume $L(S,\lambda)$ is an integrable $\gl_{\infty}$-module. Then
$L(S,\lambda)$ is in category $C_{int}$. By Lemma
\ref{lenn-eigenvalues}, $\lambda_{m}\in \N$ for all $m\in \Z$. As an
integrable $\gl_{\infty}$-module, $L(S,\lambda)$ is necessarily an
integrable $\gl_{S}$-module, containing $L(\lambda_{S})$ as a
submodule. Then it follows that $\lambda_{r}\ge \lambda_{r+1}\ge
\cdots \ge \lambda_{n}$. Thus $\lambda\in P_{+}(S)$.

Conversely, assume $\lambda\in P_{+}(S)$. By definition,
$\lambda_{m}\in \N$ for all $m\in \Z$ and $\lambda_{r}\ge
\lambda_{r+1}\ge \cdots \ge \lambda_{n}$. Notice that $E_{m,m}$ for
$m\in \Z$ are semi-simple on $L(S,\lambda)$. We now show that
$L(S,\lambda)$ is integrable, by proving that $E_{j,j+1}, E_{j+1,j}$
for $j\in \Z$ are locally nilpotent.  We shall freely use the
following observation: Suppose that $W$ is an irreducible
$\gl_{\infty}$-module. For any $p,q\in \Z$ with $p\ne q$, if
$(E_{p,q})^{k}w=0$ for some nonzero vector $w\in W$ and for some
positive integer $k$, then $E_{p,q}$ is locally nilpotent on the
whole space $W$. This simply follows from the fact that $\ad
E_{p,q}$ is locally nilpotent on $\gl_{\infty}$.

Let $v$ be a highest weight vector in $\gl_{S}$-module
$L(\lambda_{S})\subset L(S,\lambda)$.

(1) For $r\le j\le n-1$, we claim that $E_{j,j+1}$ and $E_{j+1,j}$
are locally nilpotent on $L(S,\lambda)$. With the assumption on
$\lambda$, we know that $L(\lambda_{S})$ is an (irreducible)
integrable $\gl_{S}$-module, so that $E_{j,j+1}$ and $E_{j+1,j}$ are
locally nilpotent on $L(\lambda_{S})$. Then it follows from the
simple observation.

(2) Let $p,q\in \Z$ with $p\ne q,\ q\notin S$. As $E_{p,q}v=0$, it
follows that $E_{p,q}$ is locally nilpotent on $L(S,\lambda)$.

(3) We claim that $E_{n+1,n}$ is locally nilpotent on
$L(S,\lambda)$.

 Recall that $$E_{p,q}v=0\ \ \ \mbox{
for all }p,q\in \Z\ \mbox{ with }q\notin S=\{r,r+1,\dots,n\},$$
$$E_{r,r+1}v=E_{r+1,r+2}v=\cdots =E_{n-1,n}v=0.$$ With $n+1\notin S$
we also have $E_{n,n+1}v=0$. Thus $v$ is also an $\bar{S}$-singular
vector with $\bar{S}=S\cup \{ n+1\}$. By Lemma \ref{luniversal}, we
have $L(S,\lambda)\simeq L(\bar{S},\lambda)$ which contains
$L(\lambda_{\bar{S}})$ as a $\gl_{\bar{S}}$-submodule. Notice that
$$\lambda_{i}\in \N\ \mbox{ for }i\in \Z\ \mbox{ and }\ \lambda_{r}\ge \lambda_{2}\ge \cdots \ge
\lambda_{n}\ge \lambda_{n+1}=0.$$ Then $L(\lambda_{\bar{S}})$ is an
integrable $\gl_{\bar{S}}$-module. In particular, $E_{n+1,n}$ is
locally nilpotent on $L(\lambda_{\bar{S}})$. Then it follows that
$E_{n+1,n}$ is locally nilpotent on $L(S,\lambda)$.

(4) We claim that $E_{r-1,r}$ is locally nilpotent on
$L(S,\lambda)$.

Recall that $L(\lambda_{S})$ is an integrable $\gl_{S}$-module.
Consequently, $L(\lambda_{S})$ is finite-dimensional. Let $v_{*}$ be
a lowest weight vector in $\gl_{S}$-module $L(\lambda_{S})$, so that
$$E_{j+1,j}v_{*}=0\ \ \ \mbox{ for }r\le j\le n-1.$$
We also have
$$E_{p,q}v_{*}=0\ \ \ \mbox{ for }p,q\in \Z\ \mbox{ with }\ q\notin
S.$$ In particular, we have $E_{r,r-1}v_{*}=0$. Set $\tilde{S}=S\cup
\{ r-1\}$. Then  $v_{*}$ is a lowest-weight singular vector in the
$\gl_{\tilde{S}}$-module $U(\gl_{\tilde{S}})v_{*}$. We see that
$L(S,\lambda)=U(\gl_{\infty})v_{*}$ can be naturally $\Z$-graded by
$I_{\tilde{S}}$ with $U(\gl_{\tilde{S}})v_{*}$ as the highest degree
subspace. As $L(S,\lambda)$ is irreducible, it follows that
$U(\gl_{\tilde{S}})v_{*}$ is an irreducible
$\gl_{\tilde{S}}$-module.

Let $\mu$ be the $H$-weight of $v_{*}$. By Remark
\ref{rsimple-fact}, we have $\mu_{r}\in \N$. Now consider the vector
$(E_{r-1,r})^{\mu_{r}+1}v_{*}$. As
$$E_{r,r-1}v_{*}=0\ \ \mbox{ and  }\ \
(E_{r,r}-E_{r-1,r-1})v_{*}=E_{r,r}v_{*}=\mu_{r}v_{*},$$ we have
$$E_{r,r-1}\cdot
(E_{r-1,r})^{\mu_{r}+1}v_{*}=0.$$ For $r\le j\le n-1$, since
$[E_{j+1,j},E_{r-1,r}]=0$ and $E_{j+1,j}v_{*}=0$, we have
$$E_{j+1,j}\cdot (E_{r-1,r})^{\mu_{r}+1}v_{*}=0.$$
 Furthermore, for $p,q\in \Z$ with
$q\notin S$ and $q\ne r-1$, we have
$$E_{p,q}E_{r-1,r}=E_{r-1,r}E_{p,q}-\delta_{r,p}E_{r-1,q}.$$
By induction we get $E_{p,q}\cdot (E_{r-1,r})^{\mu_{r}+1}v_{*}=0$.
Thus, $(E_{r-1,r})^{\mu_{r}+1}v_{*}$, if not zero, is another
lowest-weight singular vector in the $\gl_{\tilde{S}}$-module
$U(\gl_{\tilde{S}})v_{*}$. As $v_{*}$ and
$(E_{r-1,r})^{\lambda_{r}+1}v_{*}$ have different $H$-weights,
$(E_{r-1,r})^{\lambda_{r}+1}v_{*}$ must be zero. It then follows
that $E_{r-1,r}$ is locally nilpotent on $L(S,\lambda)$.

To summarize, we have proved that $E_{i,i+1}$ and $E_{i+1,i}$  for
$i\in \Z$ are locally nilpotent on $L(S,\lambda)$. It was known that
$H$ is semi-simple on $L(S,\lambda)$. Therefore, $L(S,\lambda)$ is
an integrable $\gl_{\infty}$-module.
\end{proof}

With Proposition \ref{pintegrability} and Theorem \ref{tmain-D},
using a standard argument (see \cite{kac1}) we obtain:

\bc{cirreducibility} Let $S=\{r,r+1,\dots,n\}$ with $r<n$ and let
$\lambda\in  P_{+}(S)$.  Then the maximal submodule of
$M(S,\lambda)$ is generated by  $E_{n+1,n}^{\lambda_{n}+1}v$,
$E_{r-1,r}^{\lambda_{r}+1}v$, and
$E_{i+1,i}^{\lambda_{i}-\lambda_{i+1}+1}v$ for $r\le i\le n-1$,
where $v$ is a highest weight vector of weight $\lambda$. \ec

To summarize we have:

\bt{tmain-classification} Let  $S=\{ r,r+1,\dots,n\}$ where $r$ and
$s$ are integers with $r<n$ and let $\lambda\in H^{*}$ be such that
${\rm supp}(\lambda)\subset S$ and $\lambda\in P_{+}(S)$. Then
$L(S,\lambda)$ belongs to category ${\cal{C}}_{int}$. On the other
hand, every irreducible $\gl_{\infty}$-module in category
${\cal{C}}_{int}$ is isomorphic to a module of this form.\et

\begin{proof} The first assertion follows from
Lemma \ref{lM(U)} and Proposition \ref{pintegrability}. Now, let $W$
be an irreducible $\gl_{\infty}$-module in ${\cal{C}}_{int}$. By
Proposition \ref{panyC}, $W\simeq L(S,U)$ for some finite subset $S$
of $\Z$ and for some irreducible $\gl_{S}$-module $U$. As $W$ is an
integrable $\gl_{\infty}$-module,  $W$ is an integrable
$\gl_{S}$-module containing $U$ as a submodule. Then $U$ is
finite-dimensional. Let $v$ be a highest weight vector in $U$ viewed
as a $\gl_{S}$-module. Noticing that $E_{q,q}v=0$ for $q\notin S$,
we see that $v$ is an $H$-weight vector of a weight $\lambda\in
H^{*}$ with ${\rm supp}(\lambda)\subset S$. It follows that
$W=U(\gl_{\infty})v\simeq L(S,\lambda)$. Furthermore, by Proposition
\ref{pintegrability} we have $\lambda\in P_{+}(S)$.
\end{proof}

\section{Decomposition of tensor product modules in category ${\cal{C}}_{int}$}
In this section, we construct the irreducible integrable
$\gl_{\infty}$-modules $L(S,\lambda)$ by using the natural module $\C^{\infty}$, and we also determine the decomposition of tensor product modules in category ${\cal{C}}_{int}$.

Let $A$ denote the polynomial algebra $\C[x_{m}\ |\ m\in \Z]$. Define $\deg
x_{m}=1$ for $m\in \Z$ to make $A$ a $\Z$-graded algebra
$$A=\oplus _{r\ge 0}A_{r}.$$
It was well known that Lie algebra $\gl_{\infty}$ naturally acts on $A$ with
$$E_{m,n}=x_{m}\frac{\partial}{\partial x_{n}}\ \ \ \mbox{ for
}m,n\in \Z.$$ It can be readily seen that $A$ is a
$\gl_{\infty}$-module in category ${\cal{C}}_{int}$. We see that
$A_{r}$ for $r\ge 0$ are submodules with $A_{0}=\C$ and $A_{1}\simeq
\C^{\infty}$. In fact, $A$ is isomorphic to the symmetric algebra
$S(\C^{\infty})$ with $v_{m}$ identified with $x_{m}$ for $m\in \Z$.

One can show that for each $r\ge 0$, $A_{r}$ is an irreducible
$\gl_{\infty}$-module. For each $n\in \Z$, let $\varepsilon_{n}$
be the linear functional on $H$ defined by
$$\varepsilon_{n}(E_{m,m})=\delta_{n,m}\ \ \ \mbox{ for }m\in \Z.$$
Clearly, $\varepsilon_{n}$ $(n\in \Z)$ are linearly independent. For
$i_{1},\dots,i_{r}\in \Z$, the $H$-weight of monomial
$x_{i_{1}}\cdots x_{i_{r}}$ is $\varepsilon_{i_{1}}+\cdots
+\varepsilon_{i_{r}}$. We see that every $H$-weight space of $A_{r}$
is $1$-dimensional. Then any nonzero submodule of $A_{r}$ must
contain a monomial of degree $r$. For $j_{1},\dots,j_{k}\in \Z,\
n_{1},\dots,n_{k}\in \N$ with $j_{1}<j_{2}<\cdots <j_{k}$ and
$n_{1}+\cdots +n_{k}=r$, we have
$$(E_{j_{1},j_{2}})^{n_{2}}\cdots (E_{j_{1},j_{k}})^{n_{k}}\left(x_{j_{1}}^{n_{1}}\cdots x_{j_{k}}^{n_{k}}\right)
=n_{2}!\cdots n_{k}!x_{j_{1}}^{r}.$$ For any $p,q\in \Z$ with $p\ne
q$, we have
$$(E_{p,q})^{r}x_{q}^{r}=r!x_{p}^{r}.$$
We also have
$$\frac{1}{m_{1}!}\cdots \frac{1}{m_{s}!}(E_{i_{1},t})^{m_{1}}\cdots (E_{i_{s},t})^{m_{s}}\cdot x_{t}^{r}
={r\choose m_{1}}{r-m_{1}\choose m_{2}}\cdots {m_{s}\choose
m_{s}}x_{i_{1}}^{m_{1}}\cdots x_{i_{s}}^{m_{s}}$$ for
$i_{1},\dots,i_{s},t\in \Z,\ m_{1},\dots,m_{s}\in \N$ with
$i_{1}<i_{2}<\cdots <i_{s}$, $m_{1}+\cdots +m_{s}=r$, and $t\ne
i_{1},\dots,i_{s}$. It then follows that $A_{r}$ is an irreducible
$\gl_{\infty}$-module.  We see that $A_{r}\simeq L(S,r\varepsilon_{1})$
with $S=\{1\}$.

On the other hand, $\gl_{\infty}$ naturally acts on the exterior
algebra $\Lambda(\C^{\infty})$, which is also an $\N$-graded module
$$\Lambda(\C^{\infty})=\bigoplus_{n\in
\N}\Lambda^{n}.$$ We have $\Lambda^{0}=\C$ and
$\Lambda^{1}=\C^{\infty}$. For $n\ge 2$, the submodule $\Lambda^{n}$
has a basis consisting of vectors
$$v_{i_{1}}\wedge\cdots \wedge v_{i_{n}}$$
for $i_{1},\dots,i_{n}\in \Z$ with $i_{1}<i_{2}<\cdots <i_{n}$. One
sees that the $H$-weight of vector $v_{i_{1}}\wedge\cdots \wedge
v_{i_{n}}$ is $\varepsilon_{i_{1}}+\dots +\varepsilon_{i_{n}}$ of
multiplicity one.  Similarly, one can show that for every $n\ge 0$,
$\Lambda^{n}$ is an irreducible $\gl_{\infty}$-module.
 For $n\ge 1$, we have $\Lambda^{n}\simeq
L(S,\varepsilon_{1}+\dots+\varepsilon_{n})$ with $S=\{
1,2,\dots,n\}$.

Note that for any permutation $\sigma$ on $\Z$,  $A_{m}$ and
$\Lambda^{n}$ are $\sigma$-invariant.

\br{rexplicit-submodule} {\em  Let $n$ be a positive integer and let
$\lambda\in H^{*}$ be such that $\lambda_{m}=0$ for $m\notin
\{1,2,\dots,n\}$, $\lambda_{i}\in \N$ for $1\le i\le n$, and
$$\lambda_{1}\ge \lambda_{2}\ge \cdots \ge \lambda_{n}.$$
Note that
$$\lambda=(\lambda_{1}-\lambda_{2})\varepsilon_{1}+(\lambda_{2}-\lambda_{3})(\varepsilon_{1}+\varepsilon_{2})
+\cdots +(\lambda_{n-1}-\lambda_{n})(\varepsilon_{1}+\cdots
+\varepsilon_{n-1})+\lambda_{n}(\varepsilon_{1}+\varepsilon_{2}+\cdots
+\varepsilon_{n}).$$ For $1\le k\le n$, set
$$w_{k}=\xi_{1}\wedge \xi_{2}\wedge \cdots \wedge \xi_{k}\in
\Lambda^{k}.$$ Furthermore, set
$$w_{\lambda}=w_{1}^{\otimes
(\lambda_{1}-\lambda_{2})}\otimes (w_{2})^{\otimes
(\lambda_{2}-\lambda_{3})}\otimes\cdots \otimes (w_{n-1})^{\otimes
(\lambda_{n-1}-\lambda_{n})} \otimes (w_{n})^{\otimes
\lambda_{n}},$$ which lies in the $\gl_{\infty}$-module
$$(\Lambda^{1})^{\otimes (\lambda_{1}-\lambda_{2})}\otimes
(\Lambda^{2})^{\otimes (\lambda_{2}-\lambda_{3})}\otimes\cdots
\otimes (\Lambda^{n-1})^{\otimes (\lambda_{n-1}-\lambda_{n})}
\otimes (\Lambda^{n})^{\otimes \lambda_{n}}.$$ Set
$S=\{1,2,\dots,n\}$. It can be readily seen that $w_{\lambda}$ is an
$S$-singular vector of weight $\lambda$. Then it follows that $U(\gl_{\infty})w_{\lambda}\simeq L(S,\lambda)$.} \er

For the rest of this section, we discuss the decomposition of
tensor product modules in category ${\cal{C}}_{int}$. In view of
Theorem \ref{tmain-D}, the tensor product of any two irreducible
modules in category ${\cal{C}}_{int}$ is completely reducible. For
example, we have
$$A_{1}\otimes A_{1}=\C^{\infty}\otimes \C^{\infty}=S^{2}(\C^{\infty})\oplus \Lambda^{2}(\C^{\infty})
=A_{2}\oplus \Lambda^{2}.$$ This typical example indicates that the
tensor product of two irreducible modules in category
${\cal{C}}_{int}$ can be a finite sum of irreducible submodules.
Next, we show that indeed
this is the case. First, we establish a technical result.

\bl{llast} Let $S$ be a finite subset of $\Z$ and let $\lambda\in
P_{+}(S)$. Let $v\in L(S,\lambda)_{\lambda}$ nonzero and set
$K_{S}^{+}={\rm span}\{ E_{p,i}\ |\ p>\max(S),\; i\in S\}$, an
abelian subalgebra of $\gl_{\infty}$. Then there exists a nonzero
$H$-weight vector $v'\in U(K_{S}^{+})v\subset L(S,\lambda)$ such
that
$$E_{p,i}v'=0\ \ \ \mbox{ for all }p\in \Z,\; i\in S.$$
\el

\begin{proof}  Let $i\in S$.
Assume $E_{i,i}v=kv$ where $k$ is a nonnegative integer (by
Lemma \ref{lenn-eigenvalues}). Then
$$E_{i,i}\left(E_{p_{1},i}\cdots E_{p_{r},i}v\right)=(k-r)E_{p_{1},i}\cdots
E_{p_{r},i}v$$ for any integers $p_{1},\dots,p_{r}$ outside $S$. As
$E_{i,i}$ has only nonnegative integer eigenvalues, we have
$$E_{p_{1},i}\cdots E_{p_{r},i}v=0\ \ \mbox{ whenever }r\ge k+1.$$
Since $K_{S}^{+}$ is abelian and $S$ is finite, it follows that
there exists a nonzero $H$-weight vector $v'\in U(K_{S}^{+})v$ such
that
\begin{eqnarray}
E_{p,i}v'=0\ \ \ \mbox{ for }p>\max(S),\; i\in S.
\end{eqnarray}
Let $p$ be a fixed integer such that $p>\max(S)$ and $E_{q,p}v'=0$
for all $q\in \Z$. For $i,j\in S$, we have $E_{p,j}v'=0$,
$E_{i,p}v'=0$, and $E_{p,p}v'=0$, so that
\begin{eqnarray}
E_{i,j}v'=E_{i,p}E_{p,j}v'-E_{p,j}E_{i,p}v'+\delta_{i,j}E_{p,p}v'=0.
\end{eqnarray}
Furthermore, let $q\notin S,\; i\in S$. As $E_{i,i}v'=0$, we have
$$\ E_{i,i}(E_{q,i}v')=E_{q,i}E_{i,i}v'-E_{q,i}v'=-(E_{q,i}v').$$
Because $E_{i,i}$ has only nonnegative integer eigenvalues, we must
have $E_{q,i}v'=0$. Therefore, we have $E_{p,i}v'=0$ for all $p\in
\Z,\; i\in S$, as desired.
\end{proof}

We shall also need the following simple fact:

\bl{lafter-complete} Let $W$ be a $\gl_{\infty}$-module in category
${\cal{C}}_{int}$ and let $S$ be a finite subset of $\Z$. Suppose
that $U$ is a $\gl_{S}$-submodule of $\Omega_{S}(W)$ such that $U$ generates $W$ as a $\gl_{\infty}$-module. If
$U=\coprod_{\alpha\in I}L(\lambda^{\alpha}_{S})$ as a
$\gl_{S}$-module with $\lambda^{\alpha}\in H^{*}$ such that ${\rm
supp}(\lambda^{\alpha}) \subset S$ for $\alpha\in I$, then
$$W\simeq \coprod_{\alpha\in I} L(S,\lambda^{\alpha}).$$
 Furthermore, $\Omega_{S}(W)=U$.
\el

\begin{proof} Since $W$ is an integrable $\gl_{S}$-module, $\Omega_{S}(W)$ as a submodule
is a direct sum of finite-dimensional irreducible $\gl_{S}$-modules.
Noticing that for any $v\in \Omega_{S}(W)$, $E_{q,q}v=0$ for $q\notin
S$, we see that any singular vector in $\Omega_{S}(W)$ viewed
as a $\gl_{S}$-module is an $S$-singular vector. Suppose that $v$ is
a singular vector in $\Omega_{S}(W)$ viewed
as a $\gl_{S}$-module of $H$-weight $\lambda^{\alpha}$ with
$\alpha\in I$. By Proposition \ref{psingular-vector},
$U(\gl_{\infty})v$ is irreducible, so that
$U(\gl_{\infty})v\simeq L(S,\lambda^{\alpha})$.
It then follows that there are singular vectors $v_{\beta}$ $(\beta\in J)$ in $U$ $(\subset \Omega_{S}(W))$
with $H$-weights $\lambda^{\beta}$ such that
$$W=\oplus_{\beta \in J}U(\gl_{\infty})v_{\beta}\simeq \oplus_{\beta \in J}L(S,\lambda^{\beta}).$$
From this, using Lemma \ref{lvacuum-irred} we get
$$\Omega_{S}(W)=\oplus_{\beta \in J}\Omega_{S}(U(\gl_{\infty})v_{\beta})=\oplus_{\beta\in J}L(\lambda^{\beta}_{S})\subset U.$$
Therefore $\Omega_{S}(W)=U$.
\end{proof}

Now, we give a decomposition into irreducible submodules of the
tensor product of any two irreducible modules in category
${\cal{C}}_{int}$.

\bt{tfinite-sum} Let $S=\{ r,r+1,\dots,n\}$ be a finite subset of
$\Z$ with $r\le n$, and let $\lambda,\mu\in P_{+}(S)$. Then
$L(S,\lambda)\otimes L(S,\mu)$ is a cyclic $\gl_{\infty}$-module.
Furthermore, there exists an integer $k$ with $k\ge n$ such that for any
integer $\bar{n}\ge k$, the
decomposition of $L(S,\lambda)\otimes L(S,\mu)$ into irreducible
$\gl_{\infty}$-submodules agrees with the decomposition of
$L(\lambda_{\bar{S}})\otimes L(\mu_{\bar{S}})$ into irreducible
$\gl_{\bar{S}}$-irreducible submodules where
$\bar{S}=\{r,r+1,\dots,\bar{n}\}$.
\et

\begin{proof} Let $u\in L(S,\lambda)_{\lambda}$ and $v\in L(S,\mu)_{\mu}$, both
nonzero. Then $u\in \Omega_{S}(L(S,\lambda)),\ v\in
\Omega_{S}(L(S,\mu))$. By Lemma \ref{llast}, there exists a nonzero
$H$-weight vector $v'\in U(K_{S}^{+})v$ such that $E_{p,q}v'=0$ for
all $p\in \Z,\; q\in S$. We now prove that $u\otimes v'$ generates
$L(S,\lambda)\otimes L(S,\mu)$ as a $\gl_{\infty}$-module. Set
$$K={\rm span}\{ E_{p,i}\ |\ p\in \Z,\ i\in S\}.$$
We have $U(K)u=U(\gl_{\infty})u=L(S,\lambda)$ (using the P-B-W
theorem) and $K\cdot v'=0$. Then
$$U(K)(u\otimes v')=U(K)u\otimes v'=L(S,\lambda)\otimes v'.$$
As $L(S,\mu)=U(\gl_{\infty})v'$, it follows that
$U(\gl_{\infty})(u\otimes v')=L(S,\lambda)\otimes L(S,\mu)$. This
proves that $L(S,\lambda)\otimes L(S,\mu)$ is cyclic on $u\otimes
v'$.

Furthermore, let $k$ be an integer larger than $n$ such that
$$v'\in \< E_{p,i}\ |\ n<p\le k,\ i\in S\>\cdot v,$$
where $\<\cdot\>$ denotes the generated subalgebra of
$U(\gl_{\infty})$. Let $\bar{n}$ be any integer larger than $k$ and set $\bar{S}=\{
r,r+1,\dots,\bar{n}\}$. Noticing that
$$u\in \Omega_{\bar{S}}(L(S,\lambda)),\ \ v'\in \Omega_{\bar{S}}(L(S,\mu)),$$
we have
\begin{eqnarray*}
L(S,\lambda)\otimes L(S,\mu)=U(\gl_{\infty})(u\otimes v')
=U(\gl_{\infty}^{(\bar{S},-)})U(\gl_{\bar{S}})(u\otimes v').
\end{eqnarray*}
As $u$ and $v$
are singular vectors in the integrable $\gl_{\bar{S}}$-modules $L(S,\lambda)$ and $L(S,\mu)$, respectively, we get
$$U(\gl_{\bar{S}})u=L(\lambda_{\bar{S}})\ \mbox{ and } \ U(\gl_{\bar{S}})v=L(\mu_{\bar{S}}).$$
Noticing that $U(\gl_{\bar{S}})v'\subset
U(\gl_{\bar{S}})v=L(\mu_{\bar{S}})$, we have
\begin{eqnarray*}
U(\gl_{\bar{S}})(u\otimes v')\subset L(\lambda_{\bar{S}})\otimes
L(\mu_{\bar{S}})\subset \Omega_{\bar{S}}\left(L(S,\lambda)\otimes
L(S,\mu)\right).
\end{eqnarray*}
By Lemma \ref{lafter-complete}  we get
\begin{eqnarray}
U(\gl_{\bar{S}})(u\otimes v')= L(\lambda_{\bar{S}})\otimes
L(\mu_{\bar{S}})=\Omega_{\bar{S}}\left(L(S,\lambda)\otimes
L(S,\mu)\right)
\end{eqnarray}
and
$$L(S,\lambda)\otimes L(S,\mu)\simeq L\left(\bar{S},L(\lambda_{\bar{S}})\otimes
L(\mu_{\bar{S}})\right).$$
This proves the second assertion.
\end{proof}

\bex{aman} {\em For an illustration, consider $A_{m}\otimes A_{n}$
with $m,n$ positive integers, where
$$A_{m}=L(S,m\varepsilon_{1})\ \ \mbox{ and }\
A_{n}=L(S,n\varepsilon_{1})$$ with $S=\{ 1\}$.  In this case, we can
show that $A_{m}\otimes A_{n}$ is cyclic on $x_{1}^{m}\otimes
x_{2}^{n}$. Set $K={\rm span}\{E_{p,1}\ |\ p\ne 1\}$. Then
$$A_{m}=U(\gl_{\infty})x_{1}^{m}=U(K)x_{1}^{m}\ \ \mbox{ and }\
\ K\cdot x_{2}^{n}=0.$$ We have
$$U(K)(x_{1}^{m}\otimes x_{2}^{n})=U(K)x_{1}^{m}\otimes x_{2}^{n}=A_{m}\otimes x_{2}^{n},$$
{}from which it follows that $U(\gl_{\infty})(x_{1}^{m}\otimes
x_{2}^{n})=A_{m}\otimes A_{n}$. Set $\bar{S}=\{1,2\}$. We have
$\gl_{\bar{S}}=\gl_{2}$, linearly spanned by $E_{i,j}$ for $1\le
i,j\le 2$, and
$$U(\gl_{\bar{S}})x_{1}^{m}=L(m\varepsilon_{1}), \ \
U(\gl_{\bar{S}})x_{2}^{n}=L(n\varepsilon_{1}).$$ For
$\gl_{2}$-modules, we have
$$L(m\varepsilon_{1})\otimes L(n\varepsilon_{1})=\bigoplus_{j=0}^{|m-n|}
L\left((m+n-j)\varepsilon_{1}+j\varepsilon_{2}\right).$$
 This gives rise to a decomposition of
$A_{m}\otimes A_{n}$ as a $\gl_{\infty}$-module $$A_{m}\otimes
A_{n}=\bigoplus_{j=0}^{|m-n|}
L\left(\bar{S},(m+n-j)\varepsilon_{1}+j\varepsilon_{2}\right).$$}
\eex


\begin{thebibliography}{DJKM2}
\bibitem[C]{ch}
V. Chari, Integrable representations of affine Lie algebras, {\em
Invent. Math.} {\bf 85} (1986) 317-335.

\bibitem[CP1]{cp1}
V. Chari and A. Pressley, New unitary representations of loop
groups, {\em Math. Ann.} {\bf 275} (1986) 87-104.

\bibitem[CP2]{cp2}
V. Chari and A. Pressley, A new family of irreducible, integrable
modules for affine Lie algebras, {\em Math. Ann.} {\bf 277} (1987)
543-562.

\bibitem[CP3]{cp3}
V. Chari and A. Pressley, Integrable representations of twisted
affine Lie algebras, {\em J. Algebra} {\bf 113} (1988) 438-464.

\bibitem[DJKM1]{djkm1}
E. Date, M. Jimbo, M. Kashiwara, T. Miwa, Operator approach to the
Kadomtsev-Petviashvili equation. Transformation groups for soliton
equations III, {\em J. Phys. Soc. Japan}, {\bf 50} (1981)
3806-3812.

\bibitem[DJKM2]{djkm2}
E. Date, M. Jimbo, M. Kashiwara, T. Miwa, A new hierarchy of soliton
equations of the KP-type. Transformation groups for soliton
equations IV, {\em Physics 4D} (1982) 343-365.

\bibitem[DKM]{dkm}
E. Date, M. Kashiwara, T. Miwa, Vertex operators and $\tau$ functions. Transformation groups for soliton
equations II, {\em Proc. Japan Acad.} {\bf 57} Ser. A (1981)
387-392.

\bibitem[FKRW]{fkrw}
E. Frenkel, V. Kac, A. Radul, W. Wang, $\mathcal{W}_{1+\infty}$ and
$\mathcal{W}(\gl_{N})$ with central charge $N$,  {\em Commun. Math.
Phys.} {\bf 170} (1995) 337-357.

\bibitem [H]{hum}
J. Humphreys, {\em Introduction to Lie Algebras and Representation
Theory}, GTM {\bf 9}, Springer, 1972.

\bibitem[JL]{jli}
C.-P. Jiang and H.-S. Li, Associating quantum vertex algebras to Lie
algebra $\widetilde{\gl}_{\infty}$,  arXiv: 1301.5833.

\bibitem[JM]{jm}
M. Jimbo and T. Miwa, Solitons and infinite diemnsional Lie algebras, {\em Publ. RIMS, Kyoyo Univ.} {\bf 19} (1983) 943-1001.

\bibitem [K]{kac1}
V. G. Kac, {\it Infinite-dimensional Lie Algebras}, 3rd ed.,
Cambridge Univ. Press, Cambridge, 1990.

\bibitem [KR]{kr}
V. Kac and A. Radul, Representation theory of the vertex algebra
$\mathcal{W}_{1+\infty}$, {\em Transf. Groups} {\bf 1} (1996) 41-70.

\bibitem[Li1]{li-mz}
H.-S. Li, On certain categories of modules for affine Lie algebras,
{\em Math. Z.} {\bf 248} (2004) 635-664.

\bibitem[Li2]{li-qva1}
H.-S. Li, Nonlocal vertex algebras generated by formal vertex
operators, {\em Selecta Math. (New Series)} {\bf 11} (2005) 349-397.

\bibitem[R]{rao1}
S. Rao, A new class of unitary representations for affine Lie
algebras, {\em J. Algebra} {\bf 120} (1989) 54-73.

\bibitem[S]{sato}
M. Sato, Soliton equations as dynamical systems on infinite
dimensional Grassmann manifolds, RIMS Kokyuroku {\bf 439} (1981)
30-46.

\bibitem[PS]{ps2000}
I. Penkov and V. Serganova, Categories of integrable $sl(\infty)$-,
$o(\infty)$-, $sp(\infty)$-modules, arXiv: 1006.2749 [math.RT].

\end{thebibliography}
\end{document}